\documentclass[journal]{IEEEtran}
\ifCLASSOPTIONcompsoc
\usepackage[nocompress]{cite}
\else
\usepackage[noadjust]{cite}
\fi

\ifCLASSINFOpdf
\usepackage[pdftex]{graphicx}
\graphicspath{{fig/}}
\usepackage[small,center]{caption}
\DeclareGraphicsExtensions{.pdf,.jpg,.png}
\else
\fi
\usepackage[cmex10]{amsmath}
\usepackage{amssymb,amsopn}
\usepackage[utf8]{inputenc}
\usepackage[english]{babel}
\newtheorem{remark}{Remark}
\newtheorem{lemma}{Lemma}

\DeclareMathOperator{\Span}{span}

\usepackage{float}
\usepackage{algorithm}
\usepackage{caption}
\usepackage{algorithmicx}
\usepackage{algpseudocode}
\def\hrf{\leavevmode\leaders\hrule height 1pt\hfill\kern0pt}
\DeclareCaptionFormat{algor}{%
\hrf\par\offinterlineskip\vskip1pt%
\textbf{#1#2}#3\offinterlineskip\hrf}
\DeclareCaptionStyle{algori}{singlelinecheck=off,format=algor,labelsep=space}
\ifCLASSOPTIONcompsoc
\usepackage[caption=false,font=normalsize,labelfon
t=sf,textfont=sf]{subfig}
\else
\usepackage[caption=false,font=footnotesize]{subfi
g}
\fi
\usepackage{alphalph}

\newtheorem{assumption}{\bf Assumption}
\begin{document}
%
\title{Aspects of 2D-Adaptive Fourier Decompositions}
%
%
%

\author{You~Gao,~Tao~Qian,~Vladimir~Temlyakov,~Long-fei~Cao
\thanks{Y. Gao and T. Qian are with Department of Mathematics, Faculty of Science and Technology, University of Macau, Macau, China. Emails: map2gao@gmail.com (Y. Gao), tqian1958@gmail.com (T. Qian)}
\thanks{V. Temlyakov is with Department of Mathematics
University of South Carolina, Columbia, SC29208, USA. Email:temlyakovusc@gmail.com}
\thanks{L. F. Cao is with Department of Mathematics, College of Sciences, China Jiliang University, Zhejiang, China. Email: feilongcao@gmail.com}
\thanks{The research was supported by Macao Government FDCT/079/2016/A2 and Multi-Year Research Grant of the University of Macau MYRG2016-00053-FST.}
\thanks{Manuscript received ......}
}

%
%

\markboth{Journal of ~,~Vol.~11, No.~4, Oct. ~2017}%
{Shell \MakeLowercase{\textit{et al.}}: for IEEE Journals}
%



\maketitle

\begin{abstract}
As a new type of series expansion, the so-called one-dimensional adaptive Fourier decomposition (AFD) and its variations (1D-AFDs) have effective applications in signal analysis and system identification. The 1D-AFDs have considerable influence to the rational approximation of one complex variable and phase retrieving problems, etc. In a recent paper, Qian developed 2D-AFDs for treating square images as the essential boundary of the 2-torus embedded into the space of two complex variables.  This paper studies the numerical aspects of multi-dimensional AFDs, and in particular 2D-AFDs, which mainly include (i) Numerical algorithms of several types of 2D-AFDs in relation to image representation; (ii) Perform experiments for the algorithms with comparisons between 5 types of image reconstruction methods in the Fourier category; and (iii) New and sharper estimations for convergence rates of orthogonal greedy algorithm and pre-orthogonal greedy algorithm. The comparison shows that the 2D-AFD methods achieve optimal results among the others. 
\end{abstract}


%
\IEEEpeerreviewmaketitle
\section{Introduction}\label{intro}
%
%
%
%
\IEEEPARstart{I}n this paper we continue to study a type of signal decomposition model, which has been developed recently and phrased as adaptive Fourier decomposition (AFD) \cite{Qian2011, Qian2010, Qian2016, Qian2014}. This type of decomposition, with a number of variations and in different contexts, may be said to be of the Fourier type, but of sparse nature. By Fourier type, we mean that the basic functions in the expansions consist of (non-tangential) boundary limits of analytic functions, and in particular, rational analytic functions. These analytic functions are defined in the two regions in, respectively, the two sides of the manifold on which the interested real-valued signals are defined. Fourier decomposition (FD) is the traditional analysis tool for expanding a signal that distributes the energy of a signal evenly in all terms. Owing to the maximal selections of the involved parameters, AFD stands as a more flexible method to adaptively represent a function at a rapid pace.

Comparing with other approximation methods, mainly including the traditional Fourier series, Fourier integral, wavelets and frame, spline and neural network, etc., the AFD ones may be said to be of the same nature as greedy algorithm (GA) (or matching pursuit). It, however, was not motivated by the GA studies \cite{Mallat1993, Davis1994, DeVore1996, Temlyakov2000, Tem1}, nor fell into any existing scope of GAs. AFD was motivated by the study of the positive analytic instantaneous frequency representation of a signal and is subtler and faster converging than any of the existing GAs in the general reproducing kernel Hilbert space context. It lately motivated a new type GA that we call pre-orthogonal greedy algorithm (Pre-OGA) \cite{Qian2016}. We present here some new results on the rate of convergence of the weak orthogonal greedy algorithm (WOGA) and the weak pre-orthogonal greedy algorithm (WPre-OGA) in a general setting -- for arbitrary Hilbert space and any dictionary. Typical results on the rate of convergence of greedy-type algorithms are as follows (see \cite{Tem1}). Under the assumption that a given element (function) $f$ belongs to a special class (usually, it is the closure of the convex hull of a symmetrized dictionary) we prove a bound for an error after $m$ iterations of the algorithm. This error bound depends on the algorithm and on the class, but does not depend on the individual element $f$ from the class. We call such error bounds {\it a priori} error bounds. A fundamentally new feature of our results is that they provide {\it a posteriori} error bound. For a given element $f$ from a special class we provide, for instance, for the WOGA error bounds, which take into account easily accessible information from realization of the previous $m$ iterations of the algorithm to give a better error estimate than the {\it a priori} one (see formula (\ref{woganeb})). 

Periodic signals may be viewed as signals defined on the unit circle. In the circle context AFD expands a given signal into a series of rational orthogonal functions, being in the span of some shifted Szeg\H{o} kernels and their derivatives in the unit disc \cite{Qian2011}. Denote by $\mathbb{D}$ the unit disc in the complex plane. For $a\in \mathbb{D}$, denote by $e_{a}$ the normalized Szeg\H{o} kernel of the unit disc $\mathbb{D}$, and $\mathcal{D}_S$ the dictionary consisting of all $e_{a}$, viz.,
$$\mathcal{D}_S=\{e_a\}_{a\in \mathbb{D}}, \quad e_{a}(z)=\frac{\sqrt{1-|a|^2}}{1-\overline{a}z}.$$
Being different from the notation with a generic dictionary context as used in section II, {\it B} \cite{Tem1}, we adopt the notation $e_a$, $a\in \mathbb{D}$, for a general dictionary element in $\mathcal{D}_S$. We note that $\mathcal{D}_S$ is redundant, and there is no orthogonality between different $e_a$'s. This notation is consistent with the literature \cite{Mi2012,Qian2016,Qian2010,Qian2013,Qian2011,Mi2016,Mo2014,Qian2014,
Mi2014,Mo2015}.
Linear combinations of $e_a$'s and their derivatives are contained in $H^2(\mathbb{D})$, the Hardy space over the unit disc $\mathbb{D}$, abbreviated as $H^2$. We aim at expanding functions into linear combinations of Szeg\H{o} kernels and their derivatives with fast convergence in energy. Fastness of the AFD algorithms is based on choosing of the parameters $a_k$ adaptively in $\mathbb{D}$ at each step through utilizing the maximal selection principle (MSP). A basic convergence rate is presented in \cite{Qian2013}. While the dictionary elements in $\mathcal{D}_S$ are not orthogonal to each other, AFD offers a clear expression of an orthonormal basis by using \textit{the generalized backward shift transform} to replace the Gram-Schmidt orthogonalization process. On one hand, the generalized backward shift transform facilitates the algorithm iterations. On the other hand, the system constructed through AFD coincides with the Takenaka-Malmquist (TM) system $\{B_k\}_{k=1}^\infty$ \cite{Takenaka1925,Walsh1935, Bultheel1999, Akccay2001}. Under the condition $\sum_{k=1}^{\infty}(1-|a_k|)=\infty$, any function $f\in H^2$ can be expanded as\begin{equation}f=\sum_{k=1}^\infty\langle f,B_k\rangle B_k,\quad\label{core}B_k(z){ }=\frac{\sqrt{1-|a_k|^2}}{1-\overline{a_k}z}\prod^{k-1}_{l=1}\frac{z-a_l}{1-\overline{a_l}z}.\end{equation}
The TM system is then complete in $H^2$ \cite{heuberger2005,Szasz1953}. We sometimes use the notation $B_k=B_{\{a_1,...,a_k\}}$ to indicate the parameter-dependence of $B_k$. The Laguerre basis \cite{szeg1939, lee1961} and two-parameter Kautz basis \cite{Broome1965} are two special cases of the TM system. The Fourier basis $\{z^k\}_{k=1}^\infty$ corresponds to the case in which all $a_k$'s in (\ref{core}) are zero. On the other hand, a sequence of the parameters $a_k$'s selected under the maximal selection principle of AFD may not satisfy the condition $\sum_{k=1}^{\infty}(1-|a_k|)=\infty$, and therefore may not define a TM basis. As compensation, the TM system, however, offers a fast converging expansion of the originally given signal. Closely related to the concept instantaneous frequency (IF) \cite{Cohen1995}, AFD, in fact, was motivated by the attempt of intrinsic positive frequency decomposition of signals. It is noted that all $B_k$ are of non-negative frequency with respect to the time variable if $a_1$ is set to be zero \cite{Qian2011,Qian2013}.

The concept of AFD is extended lately to include all sparse representations in parametrized Szeg\H{o} kernels of the context including those from the contemporary learning theory \cite{Mo2014}. For applications of AFD we refer to \cite{Mi2016,Mi2014,Mi2012,Mo2015}. We note that based on a result of M. Weiss and G. Weiss in 1962, Coifman and his colleges have been studying since 2000 a particular functional decomposition method, called unwinding Blaschke decomposition, with potential applications in sound and biomedical signal analysis \cite{Weiss1962, Nahon2000, Coifman2015a}. Their unwinding Blaschke decomposition coincides with a decomposition method independently developed in \cite{Qian2010}. Apart from the one-dimensional classical contexts, AFD has also been generalized to multi-variables with either the several complex variables setting or the Clifford algebra setting \cite{Qian2016, Qian2012, Qian2014c, Alpay2016}. It has also been generalized to matrix-valued functions \cite{Alpay2017}.

The present paper, as a compensation for the theoretical work of Qian \cite{Qian2016}, is further development and completion in several aspects of 2D-AFD including computerized algorithms and sharper estimations of Pre-OGA. We concentrate on two types of AFD in the rectangular region or equivalently on the 2-torus, viz., the 2D-AFD of the product-TM system type (Product AFD) and the two-dimensional pre-orthogonal greedy algorithm type based on the product-Szeg\H{o} dictionary (Pre-OGA, or 2D-Pre-OGA). We formulate numerical algorithms of the AFDs on the 2-torus that make the decomposition practical. In the experiments we compare the proposed algorithms with the Fourier decomposition (FD),  the greedy algorithm (GA) and the orthogonal greedy algorithm (OGA) on the product-Szeg\H{o} dictionary. We also deduce some error bound estimates of OGA and Pre-OGA in general Hilbert space, that are sharper than the previously obtained estimates. 

This paper is organized as follows. In Section II, the theory of the AFDs on the 2-torus is revised in relation to the numerical algorithms of Product AFD and Pre-OGA. In Section III we prove sharper convergence rates for OGA and Pre-OGA. In Section IV we carry out experiments on one toy and one real image data. We adopt Bhattacharyya distance, PSNR and QA (MSSIM) as numerical indicators for effectiveness measurement of the reconstruction algorithms. Discussions and conclusions are drawn in Section V.

\section{Revised 2D Signal Decomposition Methods}
In this section, we reformulate the AFD theory on the 2-torus and propose the computerized numerical realization of Product AFD and Pre-OGA. In our terminology AFD on the 2-torus includes four different types of which one is Product AFD, and the other three are of the greedy algorithm type, all being based on the product-Szeg\H{o} dictionary, namely, (standard) greedy algorithm (GA), orthogonal greedy algorithm (OGA), and pre-orthogonal greedy algorithm (Pre-OGA). For Product AFD and Pre-OGA we refer to \cite{Qian2016}, and for GA and OGA, in general,  among others, we refer to \cite{Davis1994}. The numerical realization of Product AFD and Pre-OGA are labeled as, respectively, Algorithms \ref{2dafd} and Algorithm \ref{poga}. Algorithm \ref{realsig} is prescribed as the realization of applying the above-mentioned algorithms to real-valued signals. The Fourier decomposition in the two complex variables setting (FD) is a special case of Product AFD. In fact, the realization of FD corresponds to Algorithm \ref{2dafd} with all the parameters being zero. 

Denoted by $\mathbb{T}$ the boundary of the unit disc $\mathbb{D}$ ($\mathbb{D}=\{z\in \mathbb{C}:|z|<1\}$), and by $L^2(\mathbb{T}^2)$ the Hilbert space of 2D signals with finite energy \cite{Rudin1969}. For $f,g\in L^2(\mathbb{T}^2)$, the inner product is defined $$\langle f,g\rangle =\frac{1}{4\pi^2}\int_{-\pi}^{\pi}\int_{-\pi}^{\pi}f(e^{it},e^{is}) \overline{g(e^{it},e^{is})}dtds.$$ 

All signals defined on the square $[0,2\pi)\otimes[0,2\pi)$ can be made to correspond with those defined on $\mathbb{T}^2$ under a simple change of the variables. A function $f (e^ {it}, e^ {is})$ of finite energy can be expanded as an infinite series of the basic functions in the tensor product of two Fourier systems \cite{Bracewell1986}, viz., $$f(e^{it},e^{is})=\sum_{-\infty< k,l< +\infty}c_{kl}e^{i(kt+ls)},$$ where the coefficients $c_{kl}$ are defined by 
$$c_{kl}=\frac{1}{4\pi^2}\int_{-\pi}^{\pi}\int_{-\pi}^{\pi}f(e^{it},e^{is}) \overline{e^{i(kt+ls)}}dtds \quad k,l=0, \pm1,\dots,$$
where$$\sum_{-\infty< k,l< +\infty}|c_{kl}|^2< \infty.$$
We adopt the notation $(f\otimes g)(z,w)=f(z)g(w)$.

For any real-valued signal $f$ of finite energy with its Fourier coefficients $c_{kl}, -\infty <k,l<\infty$, we define
$$f^{+,+}(e^{it},e^{is})\buildrel \text{def}\over=\sum_{k,l\geqslant 0}c_{kl}e^{i(kt+ls)},$$$$f^{+,-}(e^{it},e^{is})\buildrel \text{def}\over=\sum_{k\geqslant 0,l\leqslant0}c_{kl}e^{i(kt+ls)},$$$$f^{-,+}(e^{it},e^{is})\buildrel \text{def}\over=\sum_{k\leqslant0,l\geqslant 0}c_{kl}e^{i(kt+ls)},$$$$f^{-,-}(e^{it},e^{is})\buildrel \text{def}\over=\sum_{k\leqslant0,l\leqslant 0}c_{kl}e^{i(kt+ls)},$$
and, accordingly, the Hardy space \cite{Hardy1915}
\begin{equation}\begin{split}H^2(\mathbb{T}^2)=\{f \in L^2(\mathbb{T}^2){ }:\quad f(e^{it},e^{is})=f^{+,+}(e^{it},e^{is}), \\ \textrm{where}\; f^{+,+}(e^{it},e^{is})\buildrel \text{def}\over=\sum_{k,l\geqslant 0}c_{kl}e^{i(kt+ls)}\}.\end{split}\nonumber\end{equation}
Then for a real-valued $f\in L^2(\mathbb{T}^2)$, there holds
\begin{equation}\begin{split}\label{fc0}f(e^{it},e^{is})=&f^{+,+}(e^{it},e^{is})+f^{+,-}(e^{it},e^{is})+f^{-,+}(e^{it},e^{is})\\+&f^{-,-}(e^{it},e^{is})-F(e^{it})-G(e^{is})-c_{00}\\=&2Re\{f^{+,+}\}(e^{it},e^{is})\\+&2Re\{[f(e^{i(\cdot)},e^{-i(\cdot)})]^{+,+}\}(e^{it},e^{-is})\\ -&F(e^{it})-G(e^{is})-c_{00},\end{split}\end{equation}
where$$F(e^{it})=\frac{1}{2\pi}\int^{\pi}_{-\pi}f(e^{it},e^{is})ds,\;G(e^{is})=\frac{1}{2\pi}\int^{\pi}_{-\pi}f(e^{it},e^{is})dt.$$
To validate (\ref{fc0}) we first compute the difference
\begin{equation}\begin{split}f(e^{it},e^{is})-&[f^{+,+}(e^{it},e^{is})+f^{+,-}(e^{it},e^{is})\\+&f^{-,+}(e^{it},e^{is})+f^{-,-}(e^{it},e^{is})].\end{split}\end{equation}By considering the terms involving $k=0$ or $l=0$, we obtain\begin{equation}\begin{split}f(e^{it},e^{is})-&[f^{+,+}(e^{it},e^{is})+f^{+,-}(e^{it},e^{is})\\+&f^{-,+}(e^{it},e^{is})+f^{-,-}(e^{it},e^{is})]\\=&-\sum_{-\infty<k<\infty}c_{k0}e^{ikt}-\sum_{-\infty<l<\infty}c_{0l}e^{ils}-c_{00}\\=&-F(e^{it})-G(e^{is})-c_{00}.\nonumber\end{split}\end{equation}
Since $f$ is assumed to be real-valued, we have $$\overline{c_{kl}}=c_{-k\;-l}.$$ There holds $$\overline{f^{+,+}}=f^{-,-}, \quad\overline{f^{+,-}}=f^{-,+}.$$
Those imply \begin{equation}\begin{split}
&f^{+,+}(e^{it},e^{is})\!+f^{+,-}(e^{it},e^{is})\!+f^{-,+}(e^{it},e^{is})\!+f^{-,-}(e^{it},e^{is})\\&=2Re\{f^{+,+}\}(e^{it},e^{is})\!+\!2Re\{[f(e^{i(\cdot)},e^{-i(\cdot)})]^{+,+}\}(e^{it},e^{-is}),\nonumber\end{split}\end{equation} and thus conclude (\ref{fc0}).

We note that $F$ and $G$ are square-integrable functions on the unit circle, and the analysis of $F$ and $G$ are based on their respective Hardy space projections onto the unit circle \cite{Qian2013}. The relation in (\ref{fc0}) shows that the analysis of real-valued 2D-signals of finite energy may be reduced to that of the functions in the 1D or 2D Hardy spaces $H^2(\mathbb{D})$ and the Hardy space $H^2(\mathbb{D}^2)$. The relation (\ref{fc0}) will have crucial applications in the 2D-AFD algorithms. $f^{+,+}(e^{it},e^{is})$ and $f^{+,-}(e^{it},e^{is})$ are, in fact, the non-tangential boundary limits of the Hardy spaces in $\mathbb{D}\times \mathbb{D}$ and $\mathbb{D}\times \bar{\mathbb{D}}^c$ \cite{Garnett1981}.

\subsection{Two Dimensional Adaptive Fourier Decomposition} 
For analytic functions of several complex variables, it is not feasible to establish an AFD theory following the procedure of 1D-AFD. The critical obstacle is that there are no objects analogous with Blaschke product and backward shift operator in higher dimensions \cite{Zygmund2002}. The tensor product of two 1D-TM system (see Section \ref{intro}) is, therefore, employed. The parallel theory with the general $n$-fold tensor product of the 1D-TM systems is also available for higher dimensions. Detailed proofs can be found in \cite{Qian2016}.

Denoted by $\mathcal{B}^{\bold{a}}_N$ and $\mathcal{B}^{\bold{b}}_M$ the two finite 1D-TM systems in the unit disc with respectively the parameter vectors $\bold{a}$ and $\bold{b}$, where $\bold{a}=\{a_1,a_2,\dots,a_{N}\}$ and $\bold{b}=\{b_1,b_2,\dots,b_{M}\}$. In below, $\bold{a}$ and $\bold{b}$ are also allowed to be infinite sequences. In the infinite sequences case, we write $\mathcal{B}^{\bold{a}}$ and $\mathcal{B}^{\bold{b}}$ instead of $\mathcal{B}^{\bold{a}}_N$ and $\mathcal{B}^{\bold{b}}_M$. For $B^{\bold{a}}_k\in \mathcal{B}^{\bold{a}}_N$ and $B^{\bold{b}}_l\in \mathcal{B}^{\bold{b}}_M$, there holds $B^{\bold{a}}_k \otimes B^{\bold{b}}_l\in \mathcal{B}^{\bold{a}}_N \otimes \mathcal{B}^{\bold{b}}_M$. $\mathcal{B}^{\bold{a}}_N \otimes \mathcal{B}^{\bold{b}}_M$ is an orthonormal system in $L^2(\mathbb{T}^2)$. When the sequences $\bold{a}$ and $\bold{b}$ contain infinitely many elements, $\mathcal{B}^{\bold{a}}_N \otimes \mathcal{B}^{\bold{b}}_M$ is written as $\mathcal{B}^{\bold{a}} \otimes \mathcal{B}^{\bold{b}}$. It can be shown that when $\mathcal{B}^{\bold{a}}$ and $\mathcal{B}^{\bold{b}}$ are two bases of $H^2(\mathbb{T})$, then $\mathcal{B}^{\bold{a}}_N \otimes \mathcal{B}^{\bold{b}}_M$ is a basis of $H^2(\mathbb{T}^2)$. We will in the sequel also write $B_{\{a_1,...,a_k\}}$ for $B^{\bold{a}}_k$ and $B_{\{b_1,...,b_l\}}$ for $B^{\bold{b}}_l$. 

For a 2D signal $f\in H^2(\mathbb{T}^2)$, denote the $N$-partial sum\begin{equation}\label{sn}S_N(f)=\sum_{1\leqslant k,l\leqslant N}\langle f,B^{\bold{a}}_k \otimes B^{\bold{b}}_l\rangle B^{\bold{a}}_k \otimes B^{\bold{b}}_l=\sum^N_{k=1}D_k(f),$$ $$D_k (f) =S_k (f) -S_ {k-1} (f).\end{equation} 
For a fixed $k_0\leq N$ and the previously fixed $a_1,a_2,\dots,a_{k_0-1}$ and $b_1,b_2,\dots,b_{k_0-1}$ in $\mathbb{D}$,
\begin{equation}\|D_{k_0}(f)\|^2=\sum_{\max\{k,l\}={k_0}}|\langle f,B^{\bold{\tilde a}}_k \otimes B^{\bold{\tilde b}}_l\rangle |^2,\label{Dk}\end{equation}
where $\bold{\tilde a}=\{a_1,a_2,\dots,a_{k_0-1},a\}$ and $\bold{\tilde b}=\{b_1,b_2,\dots,b_{k_0-1},b\},a\in \mathbb{D},b\in \mathbb{D}.$

The parameters $a_{k_0}$ and $b_{k_0}$ are selected to make (\ref{Dk}) attain its maximal value among all possible selections of $a$ and $b$ inside the unit disc. Then $\bold{ a}=\{a_1,a_2,\dots,a_{k_0-1},a_{k_0}\}$ and $\bold{ b}=\{b_1,b_2,\dots,b_{k_0-1},b_{k_0}\}.$ The feasibility of such selections $a_{k_0}$ and $b_{k_0}$ is proved in \cite{Qian2016}, called the 2D maximal selection principle. 

At the above the maximal selection of each pair $(a_k, b_k)$ it is shown in \cite{Qian2016} that $$\lim_{N\to\infty} S_N = f,$$ which can be alternatively written 
\begin{equation}
f=\sum_{1\leqslant k,l\leqslant \infty}\langle f,B^{\bold{a}}_k \otimes B^{\bold{b}}_l\rangle B^{\bold{a}}_k \otimes B^{\bold{b}}_l.\label{tdafd}\end{equation}
The above mentioned theory is named as the two-dimensional AFD of the product-TM system type, or Product AFD in brief. Note that $S_N (f) $ is a $N^2$-terms approximation to $f$. The parameter selection at the $N$-th step is to obtain the maximal energy of the right-hand-side of (\ref{Dk}) consisting of $2N-1$ terms.

In numerical experiments, the signal notation $\bold{F}$ represents the discrete values of a function $f(e^{it},e^{is})\in H^2(\mathbb{T})$. Then (\ref{tdafd}) corresponds to its discrete form where the inner product is computed approximately through the complex matrix computation
\begin{equation}\label{innp2}
\langle f,B^{\bold{a}}_k \otimes B^{\bold{b}}_l\rangle\approx \frac{1}{mn}\overline{\bold{B}_{k}^{\bold{a}}}^{T}\;\bold{F}\; \overline{\bold{B}_{l}^{\bold{b}}}.
\end{equation}where $\bold{B}_{k}^{\bold{a}}$ and $\bold{B}_{l}^{\bold{b}}$ are, respectively, the discrete values of $B^{\bold{a}}_k$ and $B^{\bold{b}}_l$, and $m$ and $n$ are, respectively, the number of the discretization points of $t$ and $s$.

For representing the parameters $a$ and $b$ in (\ref{Dk}), the points set $\widetilde{\mathbb D}$ is drawn on a grid mesh in $\mathbb{D}$ :\begin{equation}\begin{split}\widetilde{\mathbb D}=\{\tilde a: \;&\tilde a=x+iy, | x|^2+|y|^2<1,\;and \; x=\frac{\tilde t}{N_t}, y=\frac{\tilde s}{N_s}\\ \textrm{where}\;& \tilde t=1-N_t,2-N_t,\dots,-1,0,1,\dots,N_t-1, \\&\tilde s=1-N_s,2-N_s,\dots,-1,0,1,\dots,N_s-1, \\&N_t, N_s\in \mathbb{Z}^+.\}\end{split}\end{equation} Then $(a,b)$ is selected in $\widetilde{\mathbb D}^2$ in our numerical experiments.

Under the maximal selection principle, our algorithm aims to obtain an $\bold F_N$, the square-sum of $N^2$ terms, to approximate $\bold F$, precisely,$$ \bold F_N=\frac{1}{mn}\sum_{1\leqslant k,l\leqslant N}\overline{\bold{B}_{k}^{\bold{a}}}^{T}\;\bold{F}\; \overline{\bold{B}_{l}^{\bold{b}}}\bold{B}_{k}^{\bold{a}} (\bold{B}_{l}^{\bold{b}})^{T}=\sum^N_{k=1}\bold D_k,$$ $$\bold D_k =\bold F_k -\bold F_ {k-1}.$$ For the previously fixed $a_1,a_2,\dots a_{k_0-1}$ and $b_1,b_2,\dots b_{k_0-1}$, $(a_{k_0},b_{k_0})$ are selected to satisfy $$(a_{k_0},b_{k_0})=\arg\max_{(a, b)\in\widetilde{\mathbb D}^2}\|\bold D_{k_0}\|^2,$$ $$\|\bold D_{k_0}\|^2=\sum_{\max\{k,l\}={k_0}}|\overline{\bold{B}_{k}^{\bold{\tilde a}}}^{T}\;\bold{F}\; \overline{\bold{B}_{l}^{\bold{\tilde b}}}|^2.$$

The above discussion is summarized in Algorithm \ref{2dafd}.
\begin{center}
\captionof{algorithm}{ \label{2dafd}2D Adaptive Fourier Decomposition of the product-TM system type (Product AFD) }
\begin{algorithmic}[1]
\Function{Afd}{$\bold F, N,\widetilde{\mathbb D}^2$}
\Statex \textit{Initialization}: 
\State iteration counter $k\gets \!1$
\State TM system $ \bold B^{\bold{\tilde a}}_1,\bold B^{\bold{\tilde b}}_1$, where $\bold{\tilde a}=a,\bold{\tilde b}=b$, $ (a, b)\in \widetilde{\mathbb D}^2$
\Statex\textit{Selection:}
\For {$k\leqslant N$ }
\State $\bold V_k\gets \|\bold D_k(\bold F)\|^2$ for all $( a, b)\in\widetilde{\mathbb D}^2$
\State $( a_k, b_k)\gets \arg\underset{( a, b)\in \widetilde{\mathbb D}^2}\max\{\bold V_k\}$
\State $\bold{ a}=\{a_1,a_2,\dots,a_{k_-1},a_k\}$
\State $\bold{ b}=\{b_1,b_2,\dots,b_{k-1},b_k\}$
\EndFor
\State $\bold F_N=\frac{1}{mn}\sum_{1\leqslant k,l\leqslant N}\overline{\bold{B}_{k}^{\bold{a}}}^{T}\;\bold F\; \overline{\bold{B}_{l}^{\bold{b}}}\bold{B}_{k}^{\bold{a}} (\bold{B}_{l}^{\bold{b}})^{T}$,
\State residuals $r_N\gets \bold F-\bold F_N$
\\ \Return $\bold a, \bold b, \bold F_N, r_N$
\EndFunction
\end{algorithmic} 
\hrf
\end{center}
\begin{remark} The composing basic element $B^{\bold a}_{k}\otimes B^{\bold b}_{l}$ becomes the Fourier basis function $e^{i(kt+ls)}$ by setting all the parameters to be zero. In such case Algorithm \ref{2dafd} reduces to multiple Fourier series expansion.\end{remark}

The above $\bold F$ is complex-valued. For a real-valued function $f\in L^2(\mathbb{T})$, the decomposition of $f$ is relying on the relation (\ref{fc0}). Denoted by $p=1,2,\dots,m$ and $q=1,2,\dots,n$ the discrete value of $t\in [0,2\pi)$ and $s\in [0,2\pi)$ respectively. Denoted by $\bold{F}^{+,+}, \bold{F}^{+,-}, \bold{F}^{-,+}, \bold{F}^{-,-}$ the discrete value of $f^{+,+},f^{+,-},f^{-,+},f^{-,-}$ in (\ref{fc0}). Denoted by $\bold{F}$ the discrete value of $f$. The discrete form of (\ref{fc0}) is:
\begin{equation}\begin{split}\bold{F}=&\bold F^{+,+}+\bold F^{+,-}+\bold F^{-,+}+\bold F^{-,-}-\bold F_0-\bold G_0-\tilde c_{00} \\=&2\operatorname{Re}\{\bold F^{+,+}\}+2\operatorname{Re}\{[\bold F(\cdot,-\cdot)]^{+,+}\}-\bold F_0-\bold G_0-\tilde c_{00}\label{Fc0} \end{split}\end{equation} where\begin{equation}\begin{split} &\bold F_0=\frac{1}{m}\sum_{p=1}^m \bold{F}(p,q),\;\bold G_0=\frac{1}{n}\sum_{q=1}^n \bold{F}(p,q),\\ &\tilde c_{00}=\frac{1}{mn}\sum_{p=1}^m\sum_{q=1}^n \bold{F}(p,q).\nonumber\end{split}\end{equation} 
It allows us to implement Algorithm \ref{2dafd} to decompose $\bold F^{+,+}$ and $[\bold F(\cdot,-\cdot)]^{+,+}$. However, the values of $\bold F^{+,+}$ and $[\bold F(\cdot,-\cdot)]^{+,+}$ are unknown. Next we show that $ \bold{F}$ contains sufficient information to obtain the decomposition of $\bold F^{+,+}$ and $[\bold F(\cdot,-\cdot)]^{+,+}$.

Denoted by $\bold E_{\bold p,k}$ the discrete points of $e^{ikt}$, $t\in [0,2\pi)$ viz., $e^{ik2\pi\frac{ p}{m}}$, $p=1,2,\dots,m$. Denoted by $\bold E_{\bold q,l}$ the discrete points of $e^{ils}$, $t\in [0,2\pi)$, viz., $e^{il2\pi\frac{ q}{n}},$ $q=1,2,\dots,n$. Let $\tilde{ \bold{F}}= \bold{F}+\bold F_0+\bold G_0+\tilde c_{00}$. It define accordingly the form $\tilde{ \bold{F}}^{+,+}$ by definition, \begin{equation}\begin{split}\tilde{ \bold{F}}^{+,+}=&\frac{1}{mn}\sum_{k,l\geqslant0}\overline{\bold E_{\bold p,k}}^{T}\;\tilde{ \bold{F}}\; \overline{\bold E_{\bold q,l}}\bold E_{\bold p,k}(\bold E_{\bold q,l})^{T}.\nonumber\end{split}\end{equation}
Since $\tilde{ \bold{F}}=\bold F^{+,+}+\bold F^{+,-}+\bold F^{-,+}+\bold F^{-,-}$,
\begin{equation}\begin{split}\tilde{ \bold{F}}^{+,+}=&\frac{1}{mn}\sum_{k,l\geqslant0}\overline{\bold E_{\bold p,k}}^{T}\;\bold{F}^{+,+}\; \overline{\bold E_{\bold q,l}}\bold E_{\bold p,k}(\bold E_{\bold q,l})^{T}\\+&\frac{1}{mn}\sum_{k,l\geqslant0}\overline{\bold E_{\bold p,k}}^{T}\;\bold{F}^{+,-}\; \overline{\bold E_{\bold q,l}}\bold E_{\bold p,k}(\bold E_{\bold q,l})^{T}\\+&\frac{1}{mn}\sum_{k,l\geqslant0}\overline{\bold E_{\bold p,k}}^{T}\;\bold{F}^{-,+}\; \overline{\bold E_{\bold q,l}}\bold E_{\bold p,k}(\bold E_{\bold q,l})^{T}\\+&\frac{1}{mn}\sum_{k,l\geqslant0}\overline{\bold E_{\bold p,k}}^{T}\;\bold{F}^{-,-}\; \overline{\bold E_{\bold q,l}}\bold E_{\bold p,k}(\bold E_{\bold q,l})^{T}\nonumber\end{split}\end{equation}
By taking the definition of $\bold F^{+,+}, \bold F^{+,-}, \bold F^{-,+}, \bold F^{-,-} $ into the equation above, there holds
\begin{eqnarray}
\left\{
\begin{array}{rcl}
\overline{\bold E_{\bold p,k}}^{T}\;\bold{F}^{+,-}\; \overline{\bold E_{\bold q,l}}=&0\\\overline{\bold E_{\bold p,k}}^{T}\;\bold{F}^{-,+}\; \overline{\bold E_{\bold q,l}}=&0\\
\overline{\bold E_{\bold p,k}}^{T}\;\bold{F}^{-,-}\; \overline{\bold E_{\bold q,l}}=&0. \end{array}
\right.
\end{eqnarray}

We therefore obtain $\tilde{ \bold{F}}^{+,+}$ is equal to $\bold F^{+,+}$, viz.,
\begin{equation}\begin{split}\tilde{ \bold{F}}^{+,+}=&\frac{1}{mn}\sum_{k,l\geqslant0}\overline{\bold E_{\bold p,k}}^{T}\;\bold{F}^{++}\; \overline{\bold E_{\bold q,l}}\bold E_{\bold p,k}(\bold E_{\bold q,l})^{T}=\bold{F}^{++}.\nonumber\end{split}\end{equation}
It means that $\bold F^{+,+}$ can be directly defined by $\tilde{ \bold{F}}$, \begin{equation}\bold F^{+,+}=\frac{1}{mn}\sum_{k,l\geqslant0}\overline{\bold E_{\bold p,k}}^{T}\;\tilde{ \bold{F}}\; \overline{\bold E_{\bold q,l}}\bold E_{\bold p,k}(\bold E_{\bold q,l})^{T}.\end{equation}
Due to the relation
\begin{eqnarray}
\left\{
\begin{array}{rcl}
\overline{\bold{B}_{k}^{\bold{a}}}^{T}\;\bold{F}^{+,-}\; \overline{\bold{B}_{l}^{\bold{b}}}=&0\\
\overline{\bold{B}_{k}^{\bold{a}}}^{T}\;\bold{F}^{-,+}\;  \overline{\bold{B}_{l}^{\bold{b}}}=&0\\
\overline{\bold{B}_{k}^{\bold{a}}}^{T}\;\bold{F}^{-,-}\;  \overline{\bold{B}_{l}^{\bold{b}}}=&0. \end{array}
\right.
\end{eqnarray}
$\bold F^{+,+}$ can be approximated by \begin{equation} \tilde{\bold F}_N^{+,+}=\frac{1}{mn}\sum_{1\leqslant k,l\leqslant N}\overline{\bold{B}_{k}^{\bold{a}}}^{T}\;\tilde{\bold{F}}\; \overline{\bold{B}_{l}^{\bold{b}}}\bold{B}_{k}^{\bold{a}} (\bold{B}_{l}^{\bold{b}})^{T},\label{fpp}\end{equation}

Similarly, for approximating $[\bold F(\cdot,-\cdot)]^{+,+}$, we need $\tilde{ \bold{F}}(\cdot,-\cdot)$. As $\tilde{ \bold{F}}$ is a period function with respect to $q$, there holds $\tilde{ \bold{F}}(p,-q)=\tilde{ \bold{F}}(p,n-q)$. Similar to (\ref{fpp}) the square-sum of $N^2$ terms to approximate $[\tilde{\bold F}(p,-q)]^{+,+}$, denoted by $[\tilde{\bold F}(\cdot,-\cdot)]^{+,+}_N$, is the decomposition of $[\tilde{ \bold{F}}(p,n-q)]^{+,+}$ into $N^2$ terms.

Product AFD on real-valued signals is summarized in the following Algorithm \ref{realsig}. In fact, any decomposition method established on $H^2$ can be applied to a real-valued function in $L^2$ through implementing Algorithm \ref{realsig}.
\begin{center}
\captionof{algorithm}{ \label{realsig}Product AFD to a real-valued signal in $L^2$}
\begin{algorithmic}[1]
\algrenewcommand\algorithmicrequire{\textbf{input:}}
\algrenewcommand\algorithmicensure{\textbf{output:}}
\Require \Statex $\bold{F}$: 2D real-valued digital signal of size $m\times n$;\Statex $\widetilde{\mathbb D}^2$: The sets of the parameters; \Statex N: Approximation level.
\Ensure \Statex $\bold{a},\bold{b},\bold{c},\bold{d}$: The parameters; \Statex $\tilde{\bold F}_N$: Approximation result; \Statex $r_N$: Relative error.
\Procedure {Afd2Image}{$\bold{F}, N,\widetilde{\mathbb D}^2$}
\Statex \textit{Initialization}: 
\State Compute $\bold F_0,\bold G_0, \tilde c_{00}$
\State $\tilde{ \bold{F}}\gets\bold{F}+\bold F_0+\bold G_0+\tilde c_{00}$
\State $\tilde{ \bold{F}}(\cdot,-\cdot)\gets \tilde{ \bold{F}}(\cdot,n-\cdot)$
\Statex\textit{Decomposition:}
\State $(\tilde{\bold F}_N^{+,+},\bold{a},\bold{b})\gets \Call{Afd}{\tilde{ \bold{F}}, N, \widetilde{\mathbb D}^2}$
\State $([\tilde{\bold F}(\cdot,-\cdot)]^{+,+}_N,\bold{c},\bold{d})\gets \Call{Afd}{\tilde{ \bold{F}}(\cdot,-\cdot), N,\widetilde{\mathbb D}^2}$
\Statex \textit{Reconstruction:}
\State Reconstruct $\tilde{\bold F}_N$ by (\ref{Fc0})
\State residuals $r_N\gets \bold F-\tilde{\bold F}_N$
\EndProcedure
\end{algorithmic}
\hrf
\end{center}
\subsection{Greedy Algorithm and Orthogonal Greedy Algorithm} 
To practically implement greedy algorithm (GA), maximal selection principle can be adopted and take a weaker form, that is by letting a weakness sequence $\tau=\{t_k\}_{k=1}^\infty$, where $t_k$ between 0 and 1. Such $t_k$ amounts to making a near-maximal selection at $k$-th step. This modification of greedy algorithm is named as weak greedy algorithm (WGA) \cite{Temlyakov2000,Tem1}. Let $\mathcal{H}$ be a general Hilbert space, and $\mathcal{D}$ a dictionary of $\mathcal{H}$, viz., for every $\psi\in \mathcal{D}$, $\|\psi\|=1$, and $\overline{\Span\mathcal{D}}=\mathcal{H}$. For $f\in \mathcal{H}$, WGA provides a greedy expansion as\begin{equation}
f=\sum_{j=1}^\infty\langle f_{j}^\tau,\psi_j^\tau\rangle\psi_j^\tau,
\end{equation}where the standard remainder $$f_{j}^\tau=f_{j}^\tau-\langle f_{j}^\tau,\psi_j^\tau\rangle\psi_j^\tau,$$ and $\psi_j^\tau\in \mathcal{D}$ is selected to satisfying$$|\langle f_{j}^\tau,\psi_j^\tau\rangle|\geqslant t_j\sup_{\psi\in \mathcal{D}}|\langle f_{j}^\tau,\psi\rangle|.$$

In the experiment part of this paper, WGA works on the product-Szeg\H{o} dictionary 
\begin{equation}\label{psd}\mathcal{A}=\{e_a\otimes e_b: \quad a,b\in \mathbb{D}\},\end{equation} where $e_a$ and $e_b$ are elements in the Szeg\H{o} dictionary $\mathcal{D}_S$ (see Section \ref{intro}). It can be verified that $\mathcal{A}$ is a dictionary of the $H^2$ space on the 2-torus. Every element function of $\mathcal{A}$ is a tensor product of two 1D parametrized Szeg\H{o} kernels. For this particular case, WGA coincides with GA. For a 2D signal $f\in H^2$, let the residual $\tilde g_1=f$. One has
\begin{equation}f(z,w)=\sum_{k=1}^{n-1}\langle \tilde{g}_k,e_{a_k}\otimes e_{b_k}\rangle e_{a_k}\otimes e_{b_k}+\tilde{g}_n(z,w),\end{equation} where $\tilde g_n$ is the standard remainder. $ (a_k, b_k)$ is selected by
\begin{equation}(a_k,b_k)=\arg\;\max_{(a, b)\in {\mathbb D}^2}\{|\langle \tilde{g}_k,e_{a}\otimes e_{b}\rangle |\quad :\quad a,b\in \mathbb{D} \}.\end{equation}Finally, it holds
\begin{equation}
f=\sum_{k=1}^{\infty}\langle \tilde{g}_k,e_{a_k} \otimes e_{b_k}\rangle e_{a_k} \otimes e_{b_k}.
\end{equation}
However, the $k$-th component is the best one term approximation of $\tilde g_k$ but it is not necessarily the best approximation from the subspace spanned by the first $k$ components. Following such an idea, the algorithm named as weak orthogonal greedy algorithm (WOGA) was established in Hilbert space \cite{DeVore1996}. For $f\in \mathcal{H}$, WOGA constructs the orthogonal standard remainder $$f_{j}^{o,\tau}=f-P_{H^\tau_j}(f), \;H^\tau_j:=\Span(\psi_1^{o,\tau}, \dots, \psi_j^{o,\tau})$$ and $\psi_j^{o,\tau}\in \mathcal{D}$ is selected to satisfying$$|\langle f_{j}^{o,\tau},\psi_j^{o,\tau}\rangle|\geqslant t_j\sup_{\psi\in \mathcal{D}}|\langle f_{j}^{o,\tau},\psi\rangle|.$$ In WOGA after the selection of a $t_n$-optimal dictionary element against the $n$-th orthogonal standard remainder we perform an orthogonal projection of the $n$-th orthogonal standard remainder onto the span of all the $n$ selected dictionary elements and thus to the smallest possible $(n+1)$-th orthogonal standard remainder.  This selection-orthogonal projection process results in speeding up the approximation, as shown in the experiments (see Section IV and V). 

Under the product-Szeg\H{o} dictionary $\mathcal{A}$, WOGA coincides with OGA. Gram-Schmidt procedure is applied along with each selected $e_{a_k}\otimes e_{b_k}$. 
For $f\in H^2(\mathbb{D}^2)$, let $u_1=e_{a_1}\otimes e_{b_1}$, one has 
\begin{equation}f=\sum_{k=1}^{\infty}\frac{\langle {g}_k,u_k\rangle }{\|u_k\|^2}u_k,\quad u_k=e_{a_k}\otimes e_{b_k}\!-\!\sum_{l=0}^{k-1}\frac{\langle e_{a_k}\otimes e_{b_k},u_l\rangle }{\|u_l\|^2}u_l,\end{equation}where $g_k$ is the $k$-th orthogonal remainder.

The first step of OGA is the same as GA. For $k>1$, it holds that $\langle {g}_k,u_k\rangle =\langle \tilde{g}_k,e_{a_k}\otimes e_{b_k}\rangle$ and $0<\|u_k\|\leqslant1$. Then at every step OGA uses the best approximation on the $\Span\{e_{a_1}\otimes e_{b_1},e_{a_2}\otimes e_{b_2},\dots,e_{a_k}\otimes e_{b_k}\}$. Defining the approximation error by the norm of the difference between original function and the $k$-partial sum of the approximation, the errors bounds of GA and OGA were compared in \cite{Tem1,Temlyakov2000,DeVore1998,DeVore1996,Davis1994}.
\subsection{Pre-Orthogonal Greedy Algorithm} 
Pre-orthogonal greedy algorithm (Pre-OGA) is proposed in \cite{Qian2016} as a new type of GA. For $f\in \mathcal{H}$, weak pre-orthogonal greedy algorithm (WPre-OGA) is given by
\begin{equation}\label{new1}
f=\sum_{j=1}^{n-1}\langle f, \xi_j\rangle \xi_j+g_n,
\end{equation}
\noindent where $\left\{\xi_1,\xi_2,\dots,\xi_j\right\}$ is the Gram-Schmidt orthogonalization of the j-system $\left\{\psi_1,\psi_2,\dots,\psi_j\right\},~ j=1,2,\dots,n-1,$ and $g_n$ is the orthogonal standard remainder. For each $j$, $\psi_j$ is selected to satisfy the Pre-Orthogonal $ t $-maximal selection principle
\begin{equation}\label{preogaselect}
\left|\langle g_j, \xi_j\rangle\right|\geq t \sup\left\{\left|\langle g_j, \xi_j^\psi\rangle\right|:\psi\in \mathcal{D}\right\},\quad t \in(0,1],
\end{equation}
where $\xi_j^\psi=0$ if $\psi \in \Span\left\{\psi_1,\psi_2,\dots,\psi_{j-1}\right\}$, otherwise $\left\{\xi_1,\xi_2,\dots,\xi_{j-1},\xi_j^\psi \right\}$ is the orthogonalization of $\left\{\psi_1,\psi_2,\dots,\psi_{j-1},\psi \right\}$. In \cite{Qian2016} we proved the convergence
\begin{equation}\label{preogainnerp}
f=\sum_{j=1}^\infty\langle f, \xi_j\rangle \xi_j.\end{equation}
Under some boundary conditions, we can have $ t  =1$, and, in the case, the corresponding algorithms is phrased as Pre-OGA. Below we will provide the analysis of such case in detail. We further assume that $\mathcal{H}$ is a reproducing kernel Hilbert space \cite{Aronszajn1950,Garnett1981}, and $\mathcal{D}$ is the dictionary consisting of the normalized reproducing kernels $\psi=\kappa_a$, where $\kappa_a$ is smoothly parametrized by $a$ in some open set $\mathrm{A}$ of $\mathbb C^m$. For $y\in \mathcal{H}$, we have, in particular,  
$$\langle f,\kappa_a\rangle=N(a)f(a),$$ where $N(a)$ is the normalizing constant to make $\|\kappa_a\|=1$. It can be shown that under the below cited Assumption \ref{assum}, (\ref{preogaselect}) is able to reach its maximum with $ t =1$ through the parameter selection. Repeating selections of the parameters are allowed, that corresponds to higher-order directional derivatives of $\kappa_{a}$ with respect to a direction at approach to $a$. 
\begin{assumption}\label{assum}
For any $j$ and $a_1,a_2,\dots,a_{j-1}\in \mathrm{A}$, $\lim_{a\to \partial\mathrm{A}}|\langle g_j,\xi_j^{\kappa_a}\rangle|=0$. The actual implementation of Pre-OGA may require the same condition for the some orders of directional derivatives $\partial_{v}\kappa_{a}$.
\end{assumption}

The underlying mechanism is as follows (also see \cite{Qian2016}). Due to Assumption 1, the next optimal $\kappa_{a_j}$ would be attainable at an interior point $a$ of $\mathrm{A}$. There are two separate cases. One is that $\kappa_{a_j}$ is not in the linear span of $\kappa_{a_1},..., \kappa_{a_{j-1}}$. In the case $\kappa_{a_j}$ is an optimal selection for $\psi_j$. The second case is that $\kappa_{a_j}$ is in the linear span of $\kappa_{a_1},..., \kappa_{a_{j-1}}$. In such case, we choose a sequence $a^{(k)}\to a_j$, where for each $k$, $\kappa_{a^{(k)}}$ is not in the span of $\kappa_{a_1},..., \kappa_{a_{j-1}}$, and in the process, where $a^{(k)}\to a_j$, the sequence $\kappa_{a^{(k)}}$ gives rise to the maximum of (\ref{preogaselect}) for $ t  =1$. In the Gram-Schmidt orthogonalization process, due to the relation
\begin{equation}
\kappa_{a_j}-\sum_{l=1}^{j-1}\langle \kappa_{a_j}, \xi_l\rangle \xi_l=0,
\end{equation}
for each $k$, we have

\begin{equation}\begin{split}&\frac{\kappa_{a^{(k)}}-\sum_{l=1}^{j-1}\langle \kappa_{a^{(k)}}, \xi_l\rangle \xi_l}{\|\kappa_{a^{(k)}}-\sum_{l=1}^{j-1}\langle \kappa_{a^{(k)}}, \xi_l\rangle \xi_l \|}\\
=&\frac{(\kappa_{a^{(k)}}-\kappa_{a_j})-\sum_{l=1}^{j-1}\langle \kappa_{a^{(k)}}-\kappa_{a_j}, \xi_l\rangle \xi_l}{\|(\kappa_{a^{(k)}}-\kappa_{a_j})-\sum_{l=1}^{j-1}\langle \kappa_{a^{(k)}}-\kappa_{a_j}, \xi_l\rangle \xi_l \|}\\
=&\frac{\frac{\kappa_{a^{(k)}}-\kappa_{a_j}}{\|a^{(k)}-a_j \|}-\sum_{l=1}^{j-1}\langle \frac{\kappa_{a^{(k)}}-\kappa_{a_j}}{\|a^{(k)}-a_j \|}, \xi_l\rangle \xi_l}{\|\frac{\kappa_{a^{(k)}}-\kappa_{a_j}}{\|a^{(k)}-a_j \|}-\sum_{l=1}^{j-1}\langle \frac{\kappa_{a^{(k)}}-\kappa_{a_j}}{\|a^{(k)}-a_j \|}, \xi_l\rangle \xi_l \|}.
\end{split}
\end{equation}
\noindent Taking limit $k\to \infty$, we have\begin{equation}\label{dirder}\begin{split}
&\lim_{k\to \infty}\frac{\kappa_{a^{(k)}}-\sum_{l=1}^{j-1}\langle \kappa_{a^{(k)}}, \xi_l\rangle \xi_l}{\|\kappa_{a^{(k)}}-\sum_{l=1}^{j-1}\langle \kappa_{a^{(k)}}, \xi_l\rangle \xi_l \|}\\
=&\frac{(\partial_{v}\kappa_{a})(a_j)-\sum_{l=1}^{j-1}\langle (\partial_{v}\kappa_{a})(a_j), \xi_l\rangle \xi_l}{\|(\partial_{v}\kappa_{a})(a_j)-\sum_{l=1}^{j-1}\langle (\partial_{v}\kappa_{a})(a_j), \xi_l\rangle \xi_l \|}\end{split}
\end{equation}
being the last term of the orthonormal $j$-system $\left\{\xi_1,\xi_2,\dots,\xi_{j-1},(\partial_{v}\kappa_{a})(a_j) \right\}$, where $\partial_{v}\kappa_{a}$ denotes the directional derivative of $\kappa_a$ along the tangential direction of $a^{(k)}\to a_j$. It means that when $\kappa_{a_j}$ is in the span of the former $\kappa_{a_l},$ $l=1,... j-1$, we use its directional derivative $(\partial_v\kappa_a)(a_j)$ instead of $\kappa_{a_j}$ in the Gram-Schmidt orthogonalization process. 

We denote by $\mathcal{D}_j$ the function set consisting of all possible normalized directional derivatives of the functions in $\mathcal{D}_{j-1}$, where we let $\mathcal{D}_0=\mathcal{D}$. The completion of the dictionary $\mathcal{D}$ is denoted $\tilde{\mathcal{D}}$, defined as $\tilde{\mathcal{D}}=\cup^\infty_{j=1}\mathcal{D}_j.$

In the $H^2(\mathbb{D})$ space under the dictionary constituted by the normalized Szeg\H{o} kernels, the Assumption 1, for the dictionary elements as well as for their directional derivatives of any order is met. We note that since the $\kappa_{a_j} =\lim_{k\to \infty}\kappa_{a^{(k)}}$ can never be a non-trivial linear combination of $\kappa_{a_l},$ $l=1,...,j-1$ \cite{Qian2016}, there follows $\kappa_{a_j} =\kappa_{a_l}$, for some $l=1,...,j-1$. In such case the above orthogonal process generates
\begin{eqnarray}\label{comszegodic}&1&,\dots,z^{m_0-1},\frac{1}{1-\overline{a_1}z},\dots,\frac{1}{(1-\overline{a_1}z)^{m_1}},\dots,\\ 
&&\frac{1}{1-\overline{a_n}z},\dots,\frac{1}{(1-\overline{a_1}z)^{m_n}},\quad n=1,2,\dots,\nonumber\end{eqnarray} where the multiple $m_j$ of $a_j$, $1\leqslant j \leqslant n$, induces the $(m-1)$-th derivatives of the corresponding Szeg\H{o} kernel. The Gram-Schmidt orthogonalization of (\ref{comszegodic}) gives rise to the TM system $\{B_1,B_2,\dots,B_n\}$ parametrized directly by $\{a_1,a_2,\dots,a_n\}$. 

Finally, we arrive at an important observation that in the 1D unit disc context the Pre-OGA under the Szeg\H{o} dictionary $\mathcal{D}_S$ is identical with 1D-AFD. 

In numerical experiments, the algorithm to realize 2D-Pre-OGA on a complex-valued digital signal $\bold F$ is Algorithm \ref{poga} as given below. Here we take $\mathcal{D}_0=\mathcal{A}$, the product-Szeg\H{o} dictionary in (\ref{psd}).
\begin{center}
\captionof{algorithm}{\label{poga}Pre-OGA on the product-Szeg\H{o} dictionary (2D Pre-OGA)} 
\begin{algorithmic}[1]
\Function{Pre-OGA}{$\bold F,N,\widetilde{\mathbb D}^2$}
\Statex \textit{Initialization}: 
\State iteration counter $k\gets \!1$
\State Dictionary $\mathcal{A}$ on the parameters set $\widetilde{\mathbb D}^2$
\State The set of the selected parameters $\bold P$ is empty, $p=0$. 
\State Approximation result $\bold F_N=0$
\Statex \textit{Selection:}
\For {$k \leqslant N$ }
\State $\bold V_k\gets |\langle \bold F,\xi^{\psi_{(\tilde a,\tilde b)}}_k\rangle |$, $\xi^{\psi_{(\tilde a,\tilde b)}}_k\in \mathcal{A}_k$.
\State $ (a_k, b_k)\gets \arg\underset{(\tilde a,\tilde b)\in \widetilde{\mathbb D}^2}\max\{\bold V_k\}$
\State $\bold F_N=\bold F_N+\langle \bold F,\xi_k\rangle\xi_k$,where $\xi_k=\xi_k^{\psi_{( a_k, b_k)}} $
\If{$(a_k, b_k)\not\in\bold P$}
\State $p=p+1$, $\bold P[p]\gets(a_k, b_k) $
\Else
\State $\psi_{(a_k,b_k)}\gets \partial_v\psi_{(a_k,b_k)}$ in $\mathcal{A}$
\EndIf
\State $\mathcal{A}_{k+1}=\{\xi_{k+1}^\psi:\xi_{k+1}^\psi=GS(\psi, \!\{\xi_l\}_{l=1}^k\!),\psi \in \mathcal{A} \!\}$.\Statex \Comment GS is Gram-Schmidt orthogonalization
\EndFor
\State residuals $r_N\gets \bold F-\bold F_N$
\\ \Return $\bold a,\bold b, \bold F_N, r_N$
\EndFunction
\end{algorithmic} 
\hrf
\end{center}

Similar with the OGA case, Pre-OGA is obviously the optimal strategy at the one-step selection over all the other strategies. This however does not theoretically guarantee the fastest convergence of Pre-OGA due to incomparability over remainders generated by different strategies. The experiments in section VI and V on the other hand show that Pre-OGA indeed practically give rise to the fastest convergence among all the compared algorithms.

\section{New Error Bound Estimations for OGA and Pre-OGA}
Let $\mathcal{H}$ be a complex Hilbert space equipped with a dictionary $\mathcal{D}$. For $M>0$, define \begin{equation}\begin{split}\label{Ham}\mathcal{H}(\mathcal{D},M)=&\{f\in \mathcal{H}: \exists \;\tilde \psi_k\in \mathcal{D},k=1,2,\dots,\\f=&\sum_{k=1}^{\infty}c_k\tilde \psi_k,\sum_{k=1}^{\infty}|c_k|\leqslant M\}.\end{split}\end{equation}
By limiting $f$ in $\mathcal{H}(\mathcal{D},M)$, the convergence rate of the energy of the remainder can be estimated. 

Under the principle of WOGA, $f\in \mathcal{H}(\mathcal{D},M)$ is of the form\begin{equation}
f=\sum_{k=1}^\infty\langle f,u_k\rangle u_k,\;u_k:=\frac{\psi^{(o)}_k-P_{H_{k-1}}(\psi^{(o)}_k)}{\|\psi^{(o)}_k-P_{H_{k-1}}(\psi^{(o)}_k)\|},\end{equation}
where $P_{H_{n-1}}$ denotes the orthogonal projection operator onto $H_{n-1}:={\Span}\{\psi^{(o)}_1,\psi^{(o)}_2,\dots,\psi^{(o)}_{n-1}\}$, $f_k=f-\sum_{l=1}^{k-1}\langle f,u_l\rangle u_l$, $\psi^{(o)}_k$ is selected such that \begin{equation}\label{mspoga}|\langle f_k,\psi^{(o)}_k\rangle|\geqslant  t _k \;\underset{\psi\in\mathcal{D}}\sup |\langle f_k,\psi\rangle |,  t _k \in(0, 1].\end{equation}

Under the principle of WPre-OGA, $f\in \mathcal{H}(\mathcal{D},M)$ is of the form\begin{equation}
f=\sum_{k=1}^\infty\langle f,B_k\rangle B_k,\end{equation}
where $B_k$ is selected out from $\{B_k^\psi\}_{\psi\in \mathcal{D}}$, such that \begin{equation}\label{msppoga}|\langle f_k,B_k\rangle |\geqslant  t _k\;\sup\{|\langle f_k,B_k^\psi\rangle|:\; \psi\in\mathcal{D}\},\quad t _k\in (0,1],\end{equation}
where $\{B_1,\dots,B_{k-1},B_k^\psi\}$ is the orthonormalization of $\{\psi^{(p)}_1,\psi^{(p)}_2\dots,\psi^{(p)}_{k-1},\psi\}$ and $B_1=\psi^{(p)}_1$.

Note that for $ t _k\in (0,1)$, the selection is always possible. For $ t _k=1$, a selection satisfying (\ref{mspoga}) and (\ref{msppoga}) may be impossible.

Under our definition, the error bound of weak orthogonal greedy algorithm (WOGA) in \cite{Temlyakov2000,Tem1} is
\begin{equation}\|f_{n}\|\leqslant \frac{M}{\sqrt{1+\sum_{k=1}^{n-1} t _k^2}},\label{wogaeb}\end{equation}
where $\{ t _k\}_{k=1}^\infty,\;0<  t _k< 1$, is a given sequence called the weak parameters.

In \cite{DeVore1996} the error bound of OGA ($ t _k=1$) is \begin{equation}\label{ogaeb}
\|f_n\|\leqslant \frac{M}{\sqrt{n}}.\end{equation}

In \cite{Qian2016}, the error bound of WPre-OGA ($ t = t _n\in (0,1)$) is
\begin{equation}\|f_n\|\leqslant \frac{R_mM}{ t }\frac{1}{\sqrt{n}},\label{wpogaeb}\end{equation}where $R_m=\max\{r_1,\dots,r_n\}, r_n=\underset{k\geq 1}\sup\{r_n(\tilde \psi_k)\}, r_n(\tilde \psi_k)=\|\tilde \psi_k-\sum^{n-1}_{l=1}\langle \tilde \psi_k, B_l\rangle B_l\|$, $\tilde \psi_k$'s are referred to (\ref{Ham}).

The above estimations will be improved in the following subsection.
\subsection{New Error Bound of WOGA}
The error bounds of WOGA can be refined as follows.

For a given $f\in \mathcal{H}(\mathcal{D},M)$ and a realization $\psi^{(o)}_1,\psi^{(o)}_2,\dots$ of the WOGA, we denote \begin{equation}\label{vn}\begin{split}v_n^2&=\|\psi^{(o)}_n-P_{H_{n-1}}(\psi^{(o)}_n)\|^2\\ &=\|\psi^{(o)}_n\|^2-\|P_{H_{n-1}}(\psi^{(o)}_n)\|^2< 1.\end{split}\end{equation} The insight for the magnitude of $v_n$ is that $\underset{n\to \infty}\lim v_n=0$. Then it follows, due to the weak maximal selection of $\psi^{(o)}_n$ \begin{equation}\begin{split}\|f_{n+1}\|^2=&\|f_{n}-\langle f_{n},u_n\rangle u_n\|^2\\=&\|f_{n}\|^2-|\langle f_{n},u_n\rangle |^2\\=&\|f_{n}\|^2-v_n^{-2}|\langle f_{n},\psi^{(o)}_n\rangle |^2\\ \leqslant &\|f_{n}\|^2-v_n^{-2}( t _n\sup_{\psi\in\mathcal{D}}|\langle f_{n},\psi\rangle |)^2.\end{split}\end{equation}
As $f_n\in \mathcal{H}$, it holds
\begin{equation}\begin{split}
|\langle f_{n},f\rangle |=&|\langle f_{n},\sum_{k=1}^\infty c_k \tilde \psi_k\rangle |\\\leqslant &\sum_{k=1}^\infty |c_k|\sup_{k\geqslant 1}|\langle f_{n},\tilde \psi_k\rangle |\\\leqslant & \sum_{k=1}^\infty |c_k|\sup_{\psi\in\mathcal{D}}|\langle f_{n},\psi\rangle |\\=&M\sup_{\psi\in\mathcal{A}}|\langle f_{n},\psi\rangle |.\label{sup0}
\end{split}\end{equation}
Namely,\begin{equation}\sup_{\psi\in\mathcal{D}}|\langle f_{n},\psi\rangle |\geqslant \frac{|\langle f_{n},f\rangle |}{M}= \frac{\|f_{n}\|^2}{M}.\label{sup}\end{equation}The last equality is due to the relation of the orthogonality $\langle f_n,f-f_n\rangle=0.$
Then we have
\begin{equation}\begin{split}\|f_{n+1}\|^2\leqslant & \|f_{n}\|^2(1-(\frac{ t _n}{v_n})^2\frac{\|f_{n}\|^2}{M^2}).\end{split}\end{equation}
It remains to use the following lemma \cite{Temlyakov2000,Tem1}.
\begin{lemma}\label{amb}
Let $\{a_m\}_{m=1}^\infty$ be a sequence of non-negative numbers satisfying the inequalities, $$a_1\leqslant A, t_1=1, a_m\leqslant a_{m-1}(1-\frac{t_{m-1}^2a_{m-1}}{A}),\quad m=2,3,\dots$$
Then we have for each $m$, $$a_m\leqslant \frac{A}{1+\sum_{k=1}^{m-1}t_k^2}.$$
\end{lemma}
By using Lemma \ref{amb}, we obtain the newer error bound for WOGA \begin{equation}\|f_{n}\|\leqslant M(1+\sum_{k=1}^{n-1}(\frac{ t _k}{v_k})^2)^{-\frac{1}{2}}.\label{woganeb}\end{equation}
The last estimation is an improvement of (\ref{wogaeb}) as (\ref{vn}) holds.

The new bound (\ref{woganeb}) suggests the following additional criterion for selection of an element $\psi_m^{(o)}\in\mathcal D$ at the $m$-th iteration. At the greedy step of the $m$-th iteration choose out of those $\psi_m^{(o)}$, which satisfy the weak greediness assumption (\ref{mspoga}), the one with the smallest $v_m$. 
\subsection{New Error Bounds of WPre-OGA}
In this subsection, we will deduce a sharper error bound estimation for WPre-OGA similar to (\ref{wpogaeb}).
\subsubsection{Estimation in $\mathcal{H}(\mathcal{D},M)$}
We have the following equality:
\begin{equation}\begin{split}\label{n1}\|f_{n+1}\|^2=&\|f_{n}-\langle f_{n},B_n\rangle B_n\|^2\\=&\|f_{n}\|^2-|\langle f_{n},B_n\rangle |^2.\end{split}\end{equation} 
where $B_k=\frac{\psi_k-P_{H_{n-1}}(\psi_k)}{\|\psi_k-P_{H_{n-1}}(\psi_k)\|}.$
\\The estimation of $|\langle f_{n},B_n\rangle |$ can be given from the selection principle as
\begin{equation}\begin{split}\label{msp}|\langle f_{n},B_n\rangle |\geqslant & t _n\sup_{\psi\in \mathcal{D}}|\langle f_{n},B_n^{\psi}\rangle |\\\geqslant & t _n\sup_{k\geqslant 1}|\langle f_{n},B_n^{\tilde \psi_k}\rangle | \\ =&  t _n\sup_{k\geqslant 1}\frac{|\langle f_{n},\tilde \psi_k\rangle |}{r_n(\tilde \psi_k)} \\ \geqslant & \frac{ t _n}{r_n}\sup_{k\geqslant 1}|\langle f_{n},\tilde \psi_k\rangle | \end{split}\end{equation}
 and $r_n=\underset{k\geqslant 1}\max\{r_n(\tilde \psi_k)\}$, $r_n(\tilde \psi_k)=\|\tilde \psi_k-P_{H_{n-1}}(\tilde \psi_k)\|$.

Similar to the first inequality in (\ref{sup0}), for $f\in \mathcal{H}(\mathcal{D},M)$ we have $$\sup_{k\geqslant 1}|\langle f_{n},\tilde \psi_k\rangle |\geqslant \frac{|\langle f_{n},f\rangle |}{M}= \frac{\|f_{n}\|^2}{M}.$$
Taking the result into (\ref{n1}) and (\ref{msp}), it is obtained that \begin{equation}\begin{split}\|f_{n+1}\|^2\leqslant & \|f_{n}\|^2(1-(\frac{ t _n}{r_n})^2\frac{\|f_{n}\|^2}{M^2}).\end{split}\end{equation}
The bound of WPre-OGA under Lemma \ref{amb} is as\begin{equation}\|f_{n}\|\leqslant M(1+\sum_{k=1}^{n-1}(\frac{ t _k}{r_k})^2)^{-\frac{1}{2}},\; t _k\in(0,1].\label{wpreneb}\end{equation}

At this point, we would like to comment that WPre-OGA is able to reach $ t _n=1$ at each step, if the dictionary is suitably enlarged. We adopt the concept complete dictionary (\cite{Qian2016}) induced from $\mathcal{D}$ into consideration. The complete dictionary induced from $\mathcal{D}$ is \begin{equation}\tilde{\mathcal{D}}=\underset{k=0}{\overset{\infty}\cup}\mathcal{D}_k,\end{equation}where $\mathcal{D}_k,\;k=1,2,\dots$, is the function set consisting of all possible normalized directional derivatives of the functions in $\mathcal{D}_{k-1}$ and $\mathcal{D}_0=\mathcal{D}$. In the following content, we will use a new sub-classes of functions in $\mathcal{H}$ defined in terms of the complete dictionary $\tilde{\mathcal{D}}$, for the error bound estimation of Pre-OGA.
\subsubsection{Estimation in ${\mathcal{H}}(\tilde{\mathcal{D}},N)$}
The sub-classes of functions are given as \begin{equation}\begin{split}&{\mathcal{H}}(\tilde{\mathcal{D}},N)=\{f\in \mathcal{H},\;\exists\; \tilde B_k,\; k=1,2,\dots,\\&f=\sum_{k=0}^\infty d_k\tilde B_k, \;\tilde B_k\in \underset{l=0}{\overset{k}\cup}\mathcal{D}_l,\;\sum_{k=0}^\infty|d_k|\leqslant N\}.\end{split}\end{equation} 
With regards to the selection principle (\ref{msp}), now the supremum value is attainable at an element $\tilde \psi\in \tilde{\mathcal{D}}$. This fact needs a proof, and in fact, is proved in \cite{Qian2016}, in which the completed dictionary is introduced. And, we have,
\begin{equation}\begin{split}|\langle f_n,B_n^{\tilde \psi}\rangle |= \max_{\psi\in \mathcal{D}}|\langle f_n,B_n^{\psi}\rangle |.\end{split}\end{equation}

\begin{table*}[t]
\centering
\caption{The real part of $\bold F_1,\;\bold F_9$ and $\bold F_{25}$ and remainders after $25$-terms approximation}
\label{f5img}
\begin{tabular}{cccccc}
\hline
Algorithm&FS&GA&OGA&Product AFD&2D-Pre-OGA\\
\hline
\noalign{\vskip 1mm}
\begin{tabular}{c}$\bold F_1$\\\noalign{\vskip20mm}\end{tabular} &\includegraphics[width=1.1in]{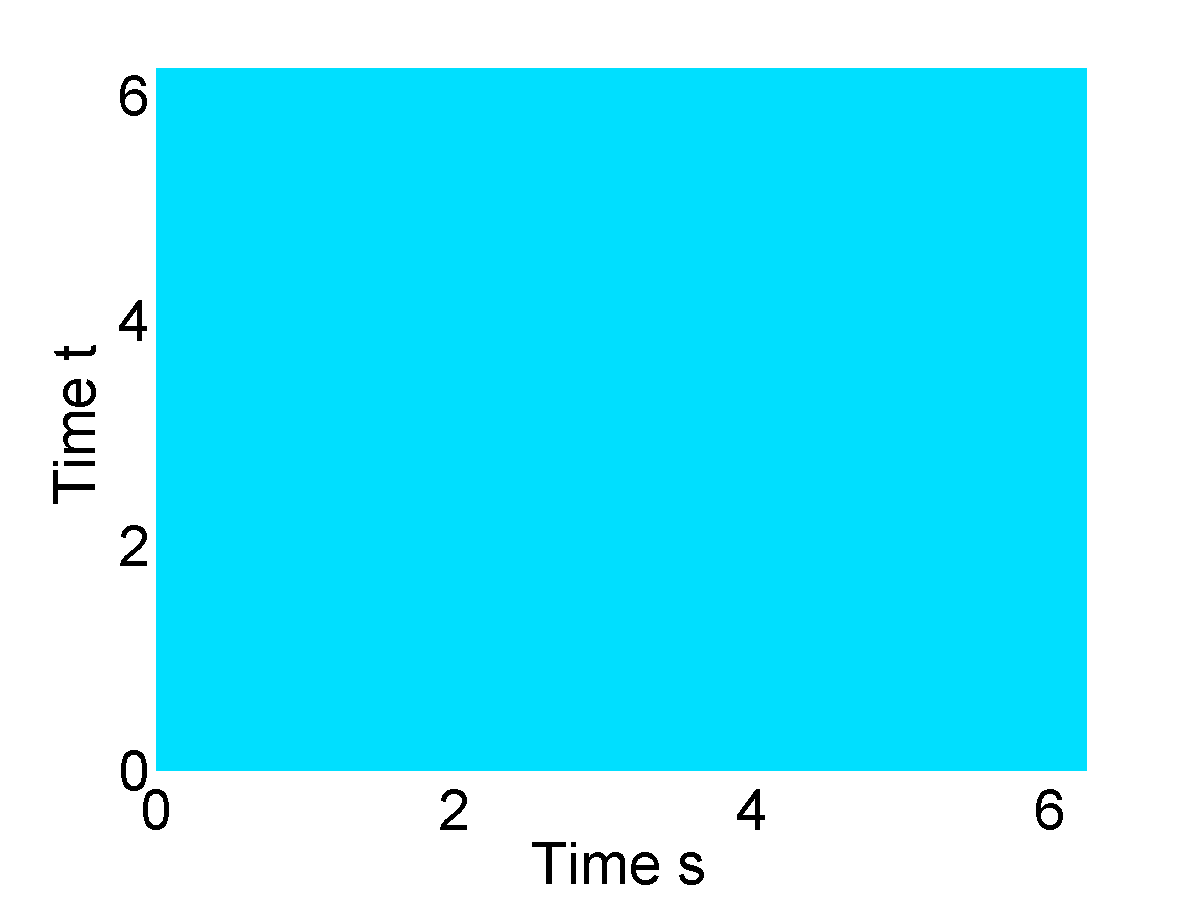}&\includegraphics[width=1.1in]{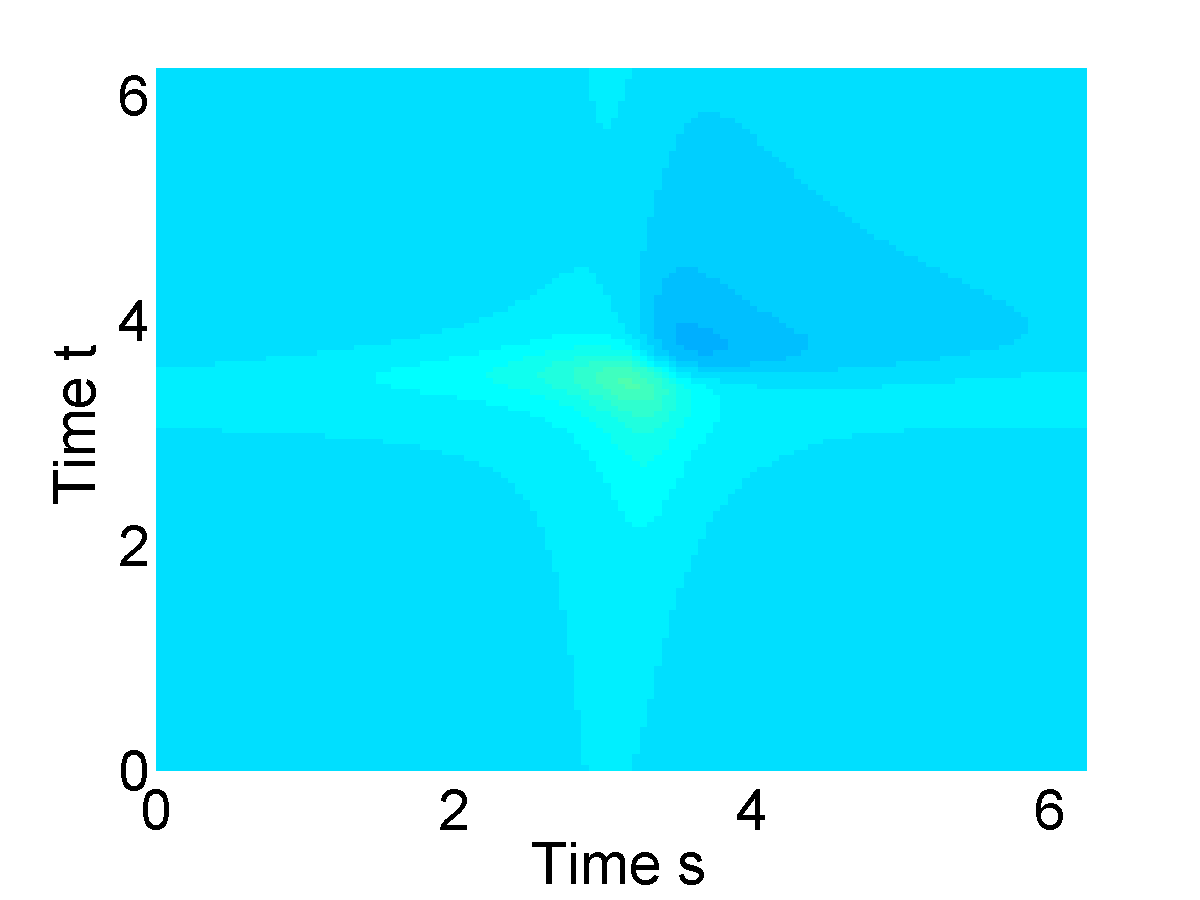}&\includegraphics[width=1.1in]{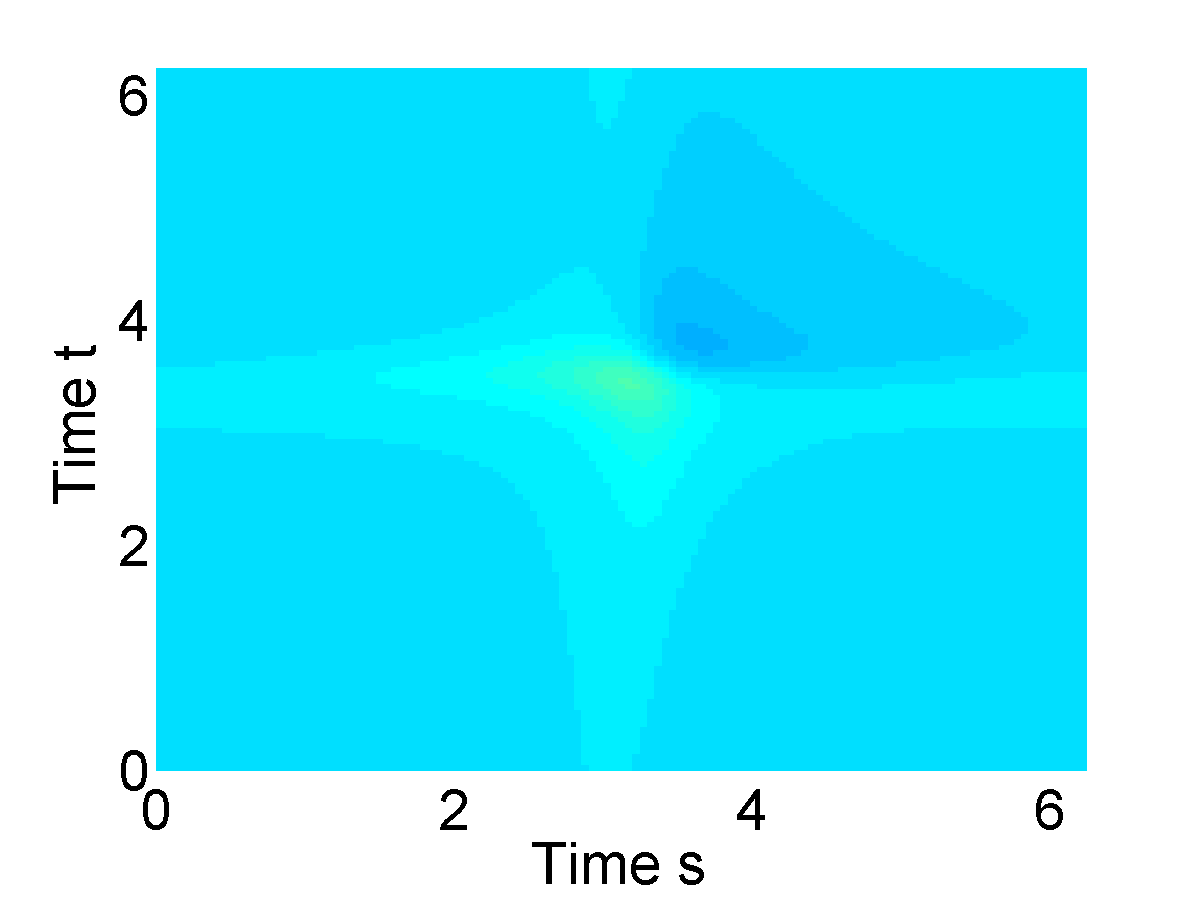}&\includegraphics[width=1.1in]{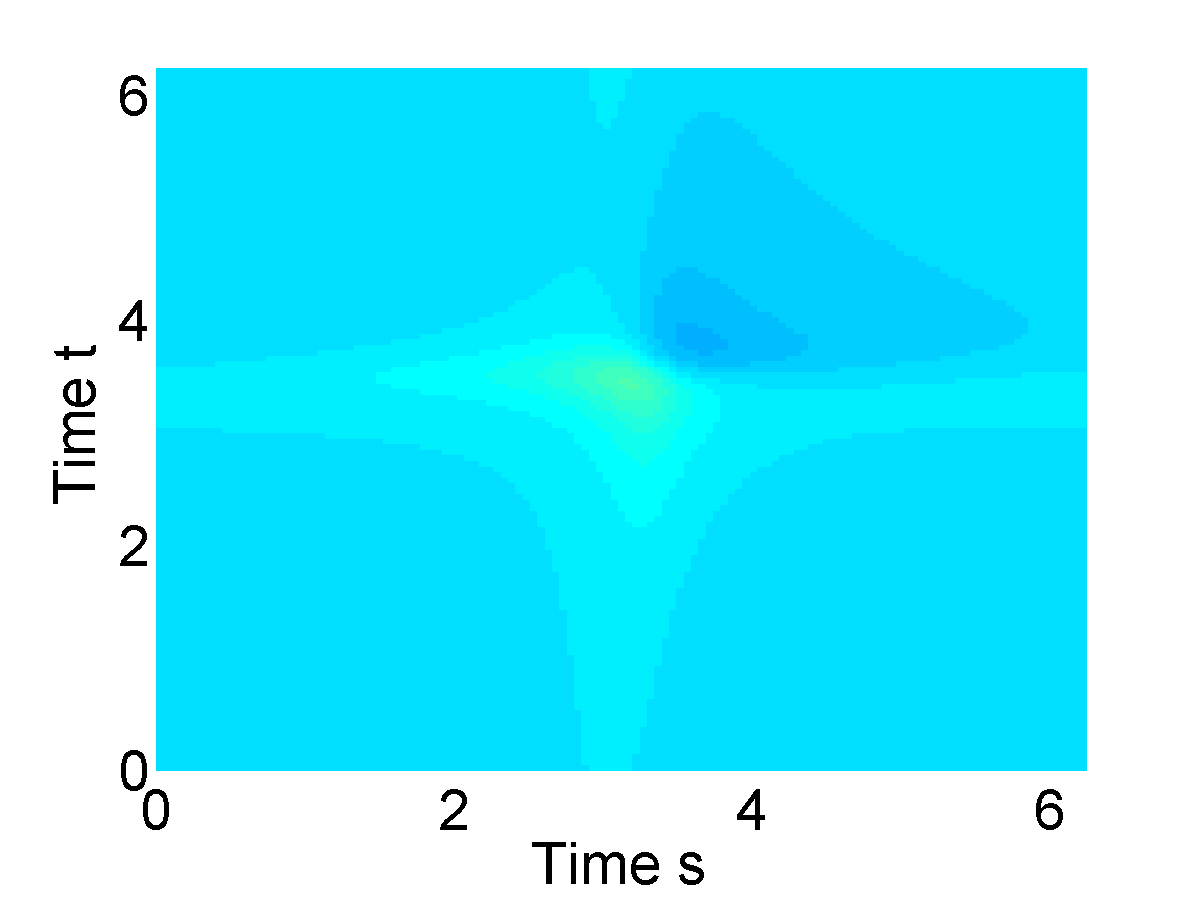}&\includegraphics[width=1.1in]{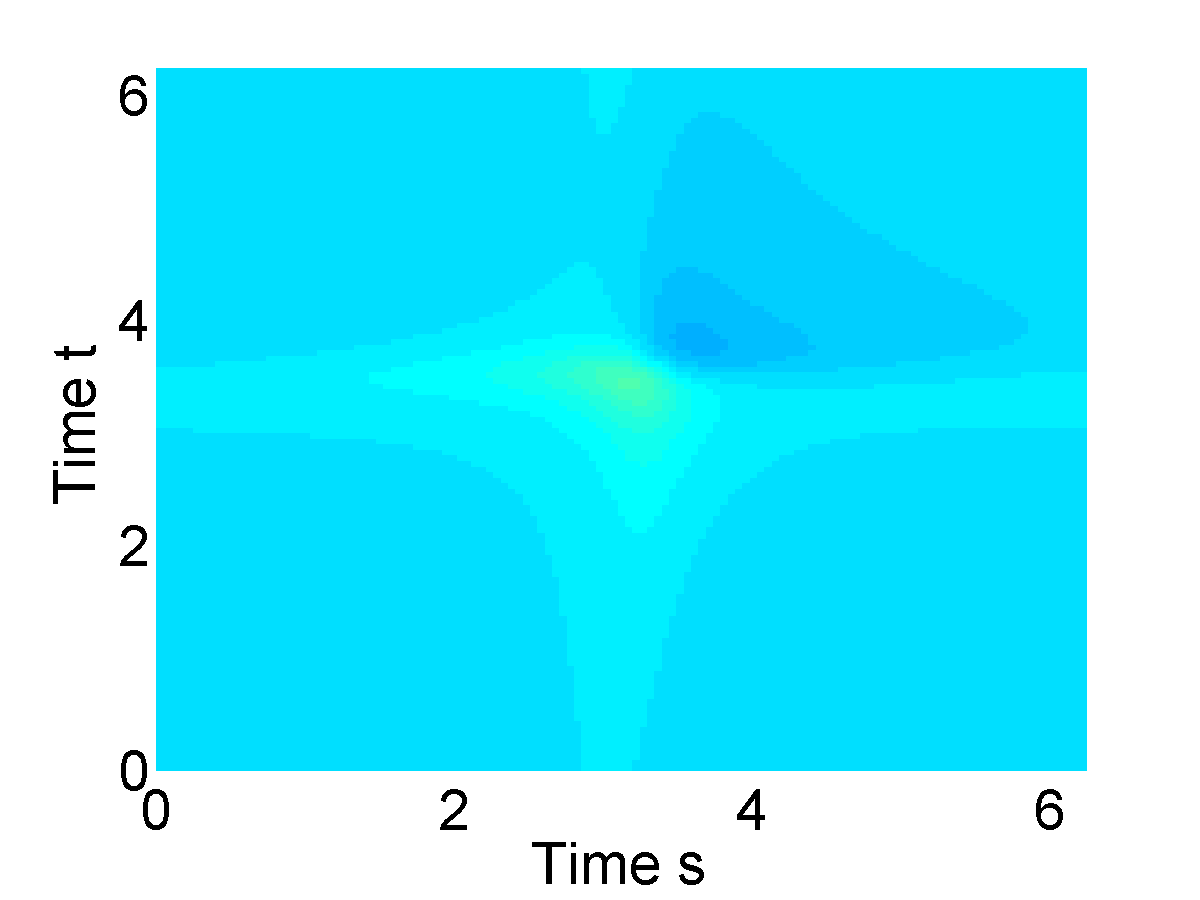}\\
\noalign{\vskip-10mm}
\hline
\noalign{\vskip 1mm}
\begin{tabular}{c}$\bold F_9$\\\noalign{\vskip20mm}\end{tabular}&\includegraphics[width=1.1in]{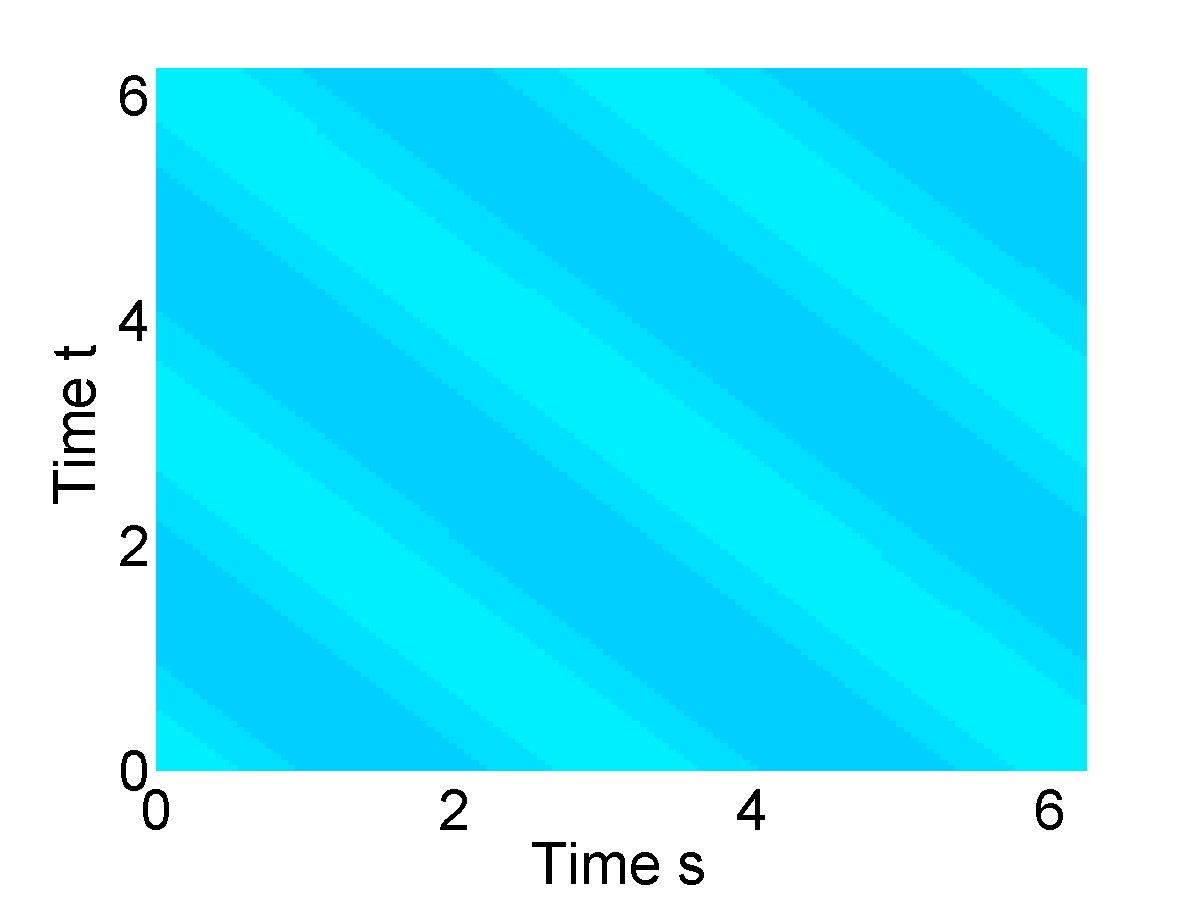}&\includegraphics[width=1.1in]{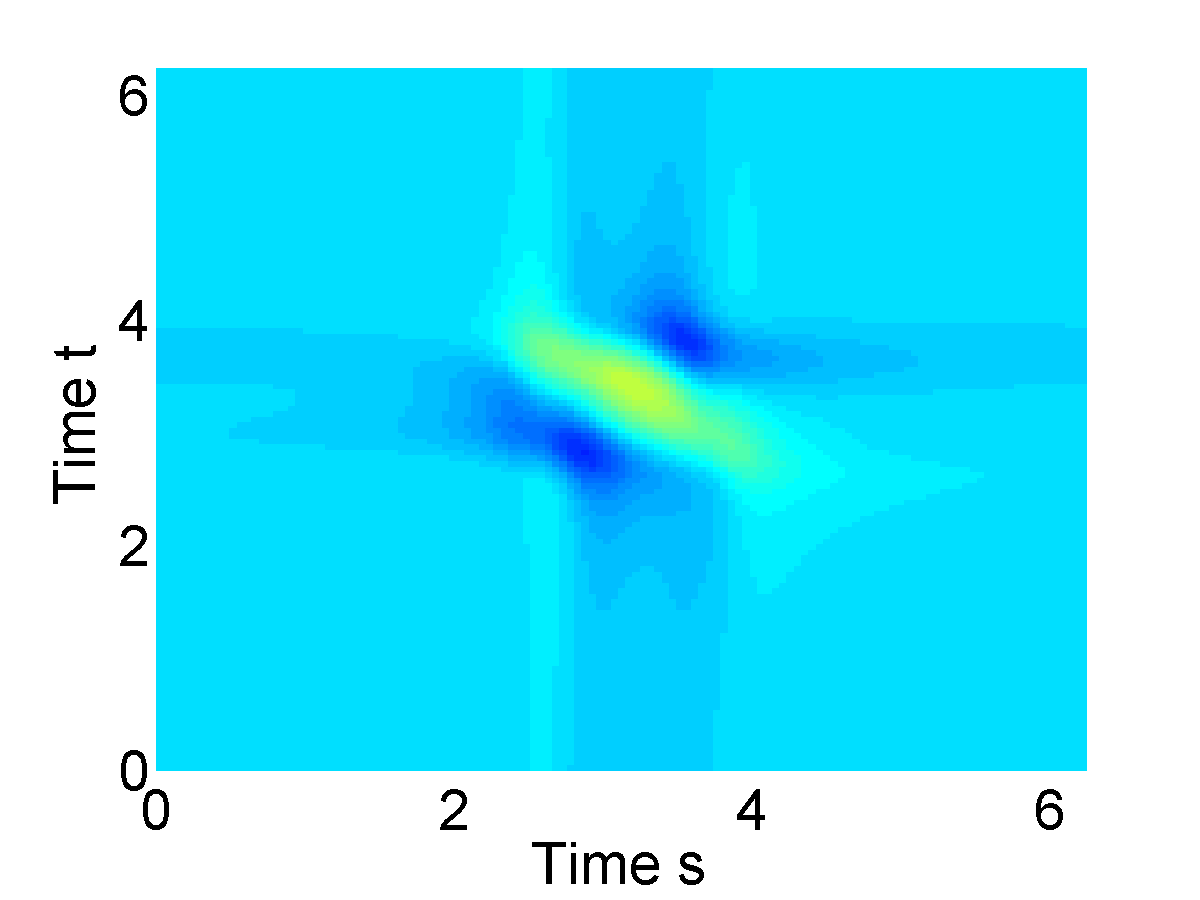}&\includegraphics[width=1.1in]{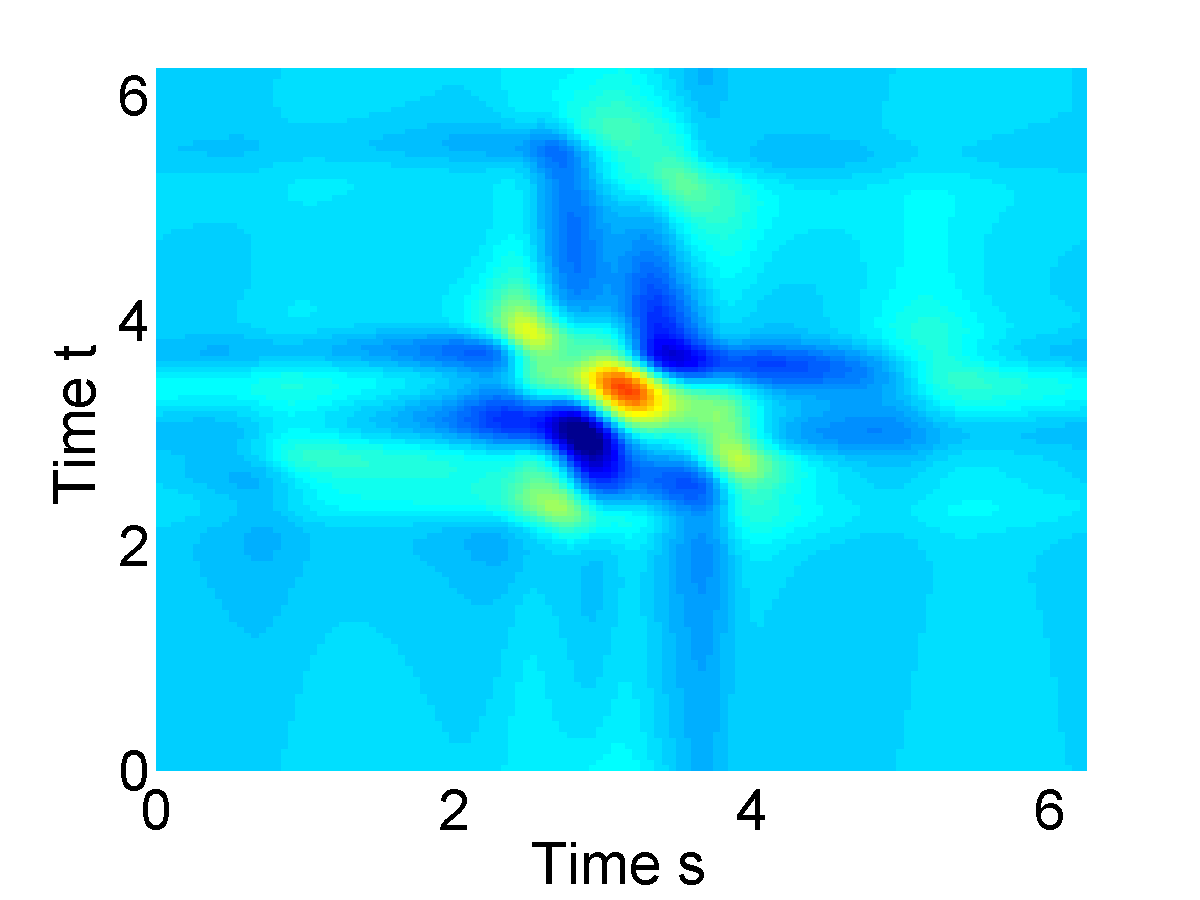}&\includegraphics[width=1.1in]{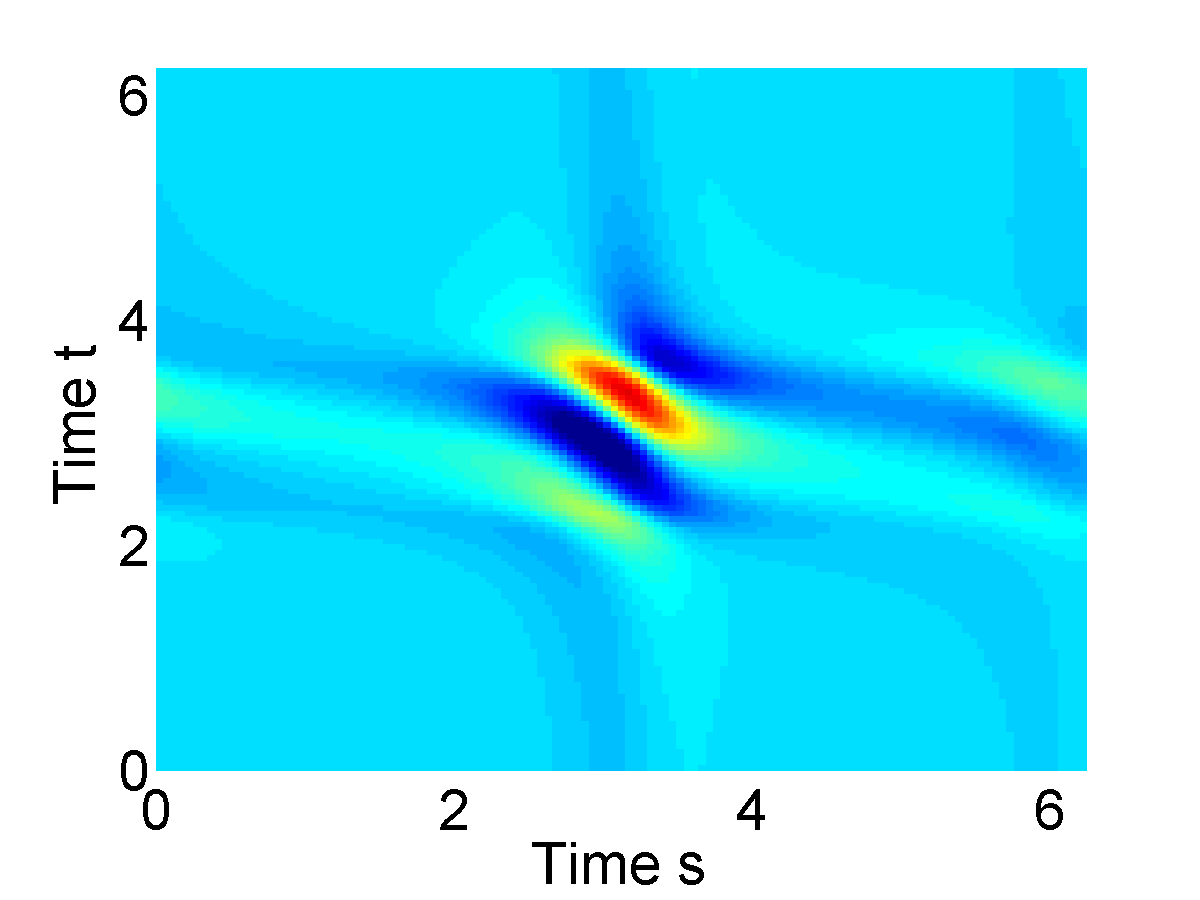}&\includegraphics[width=1.1in]{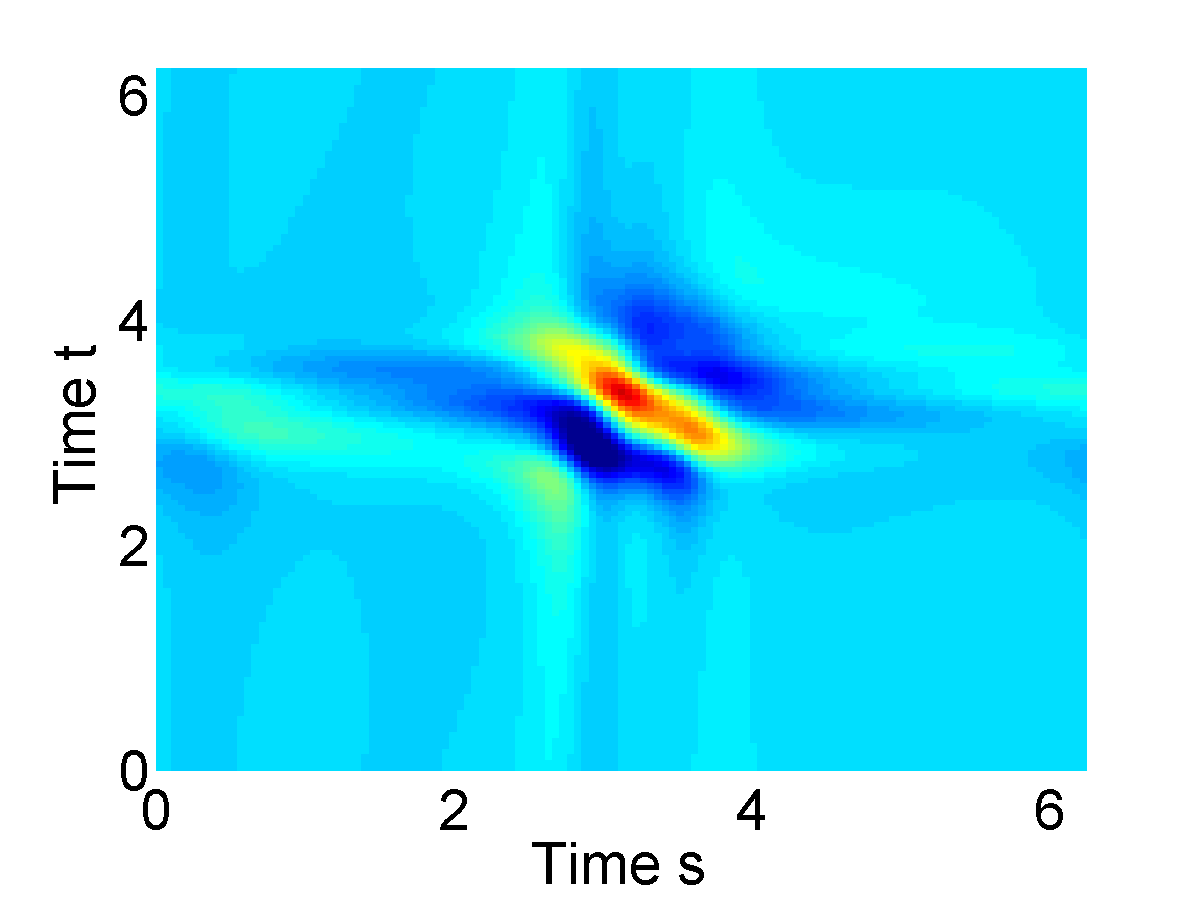}\\
\noalign{\vskip-10mm}
\hline
\noalign{\vskip 1mm}
\begin{tabular}{c}$\bold F_{25}$\\\noalign{\vskip20mm}\end{tabular}&\includegraphics[width=1.1in]{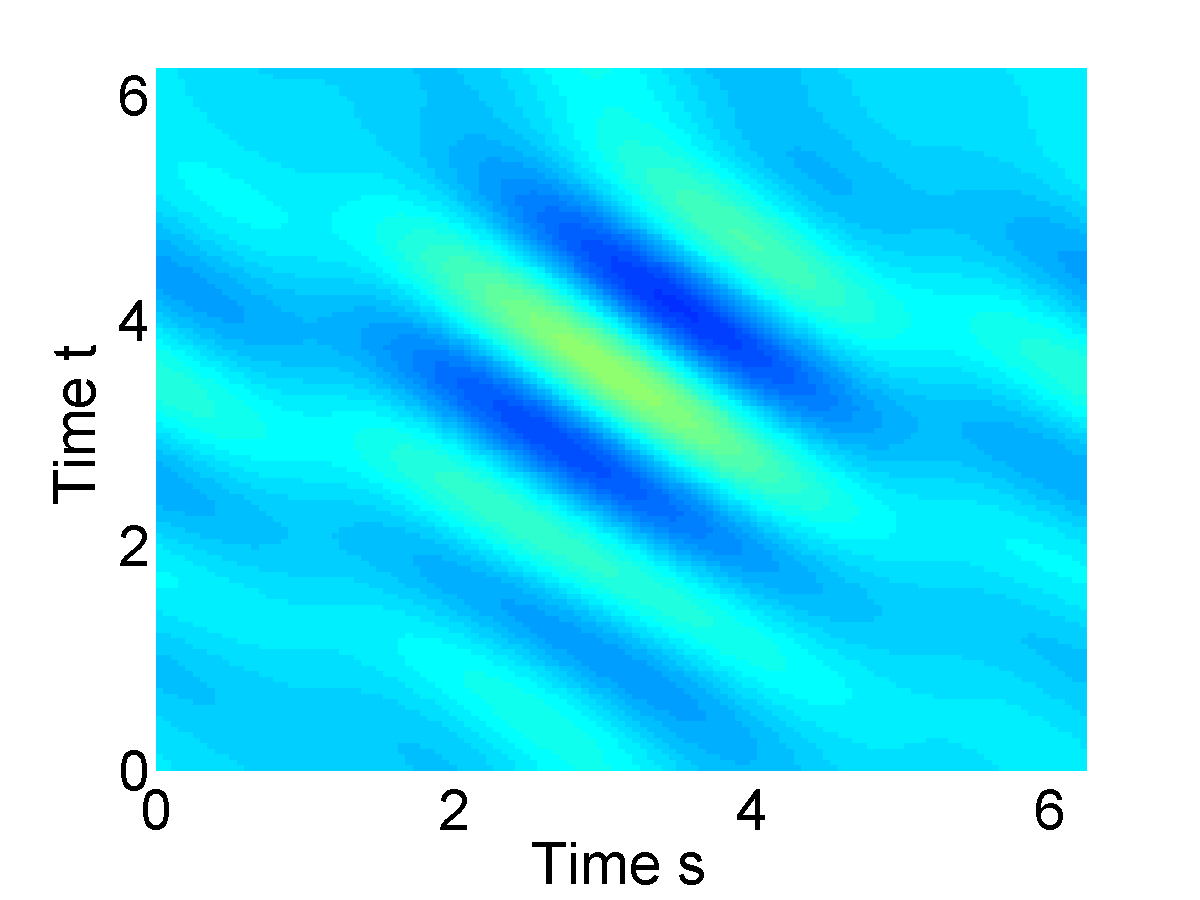}&\includegraphics[width=1.1in]{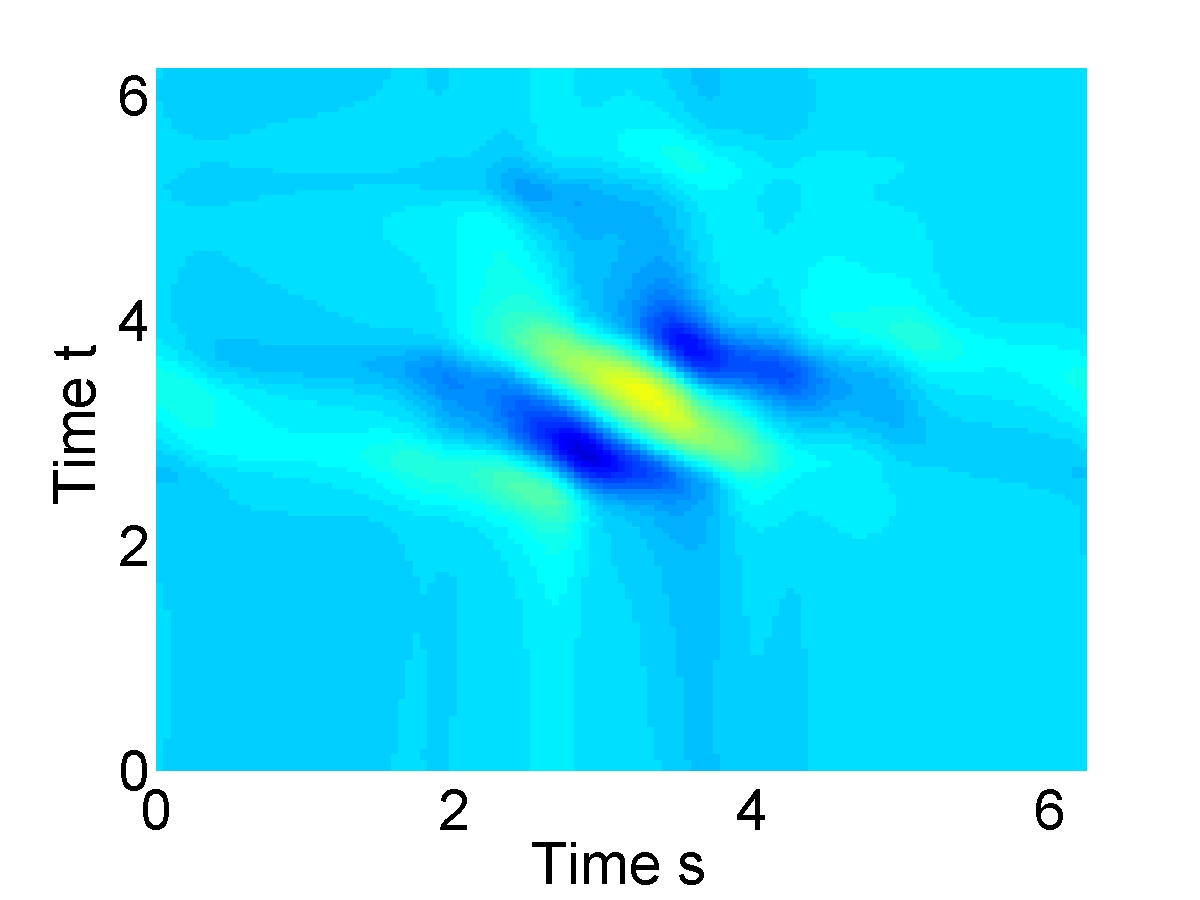}&\includegraphics[width=1.1in]{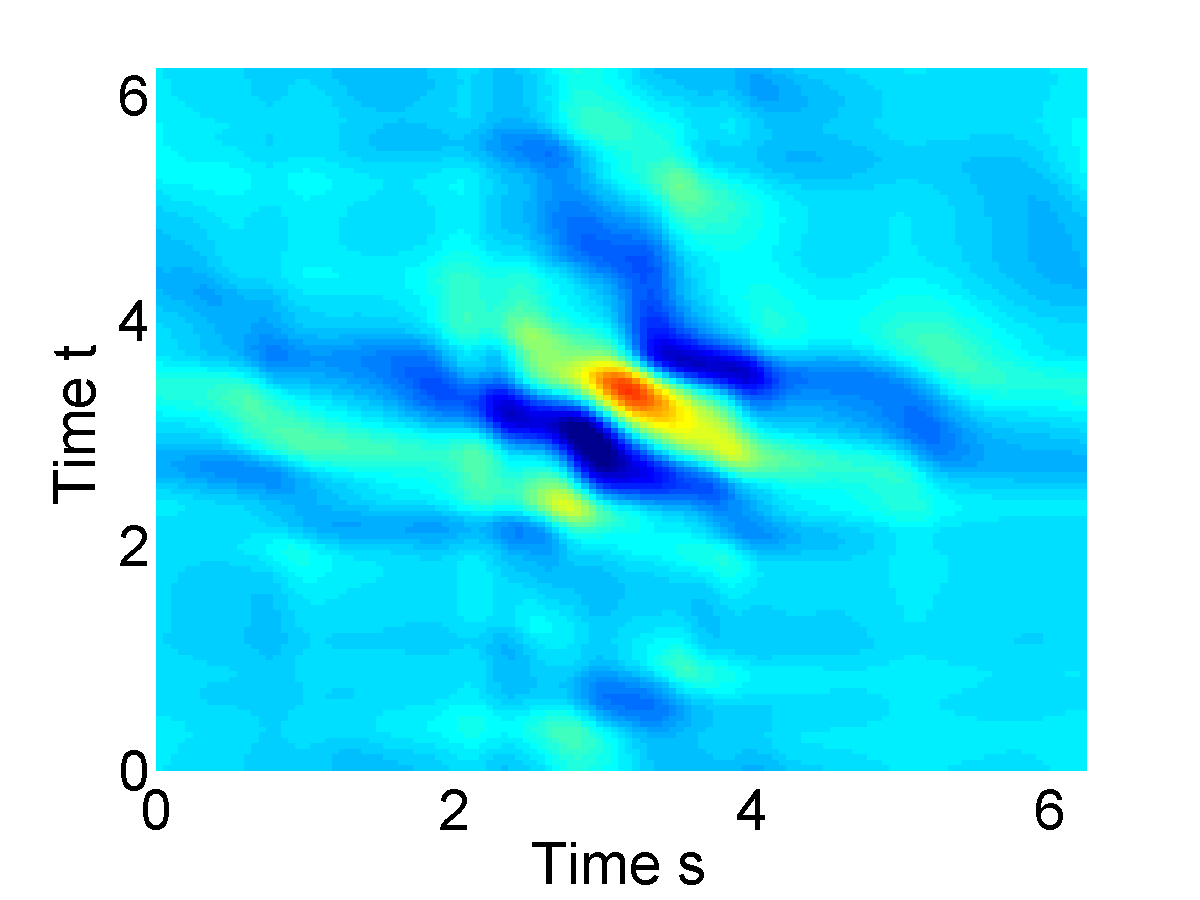}&\includegraphics[width=1in]{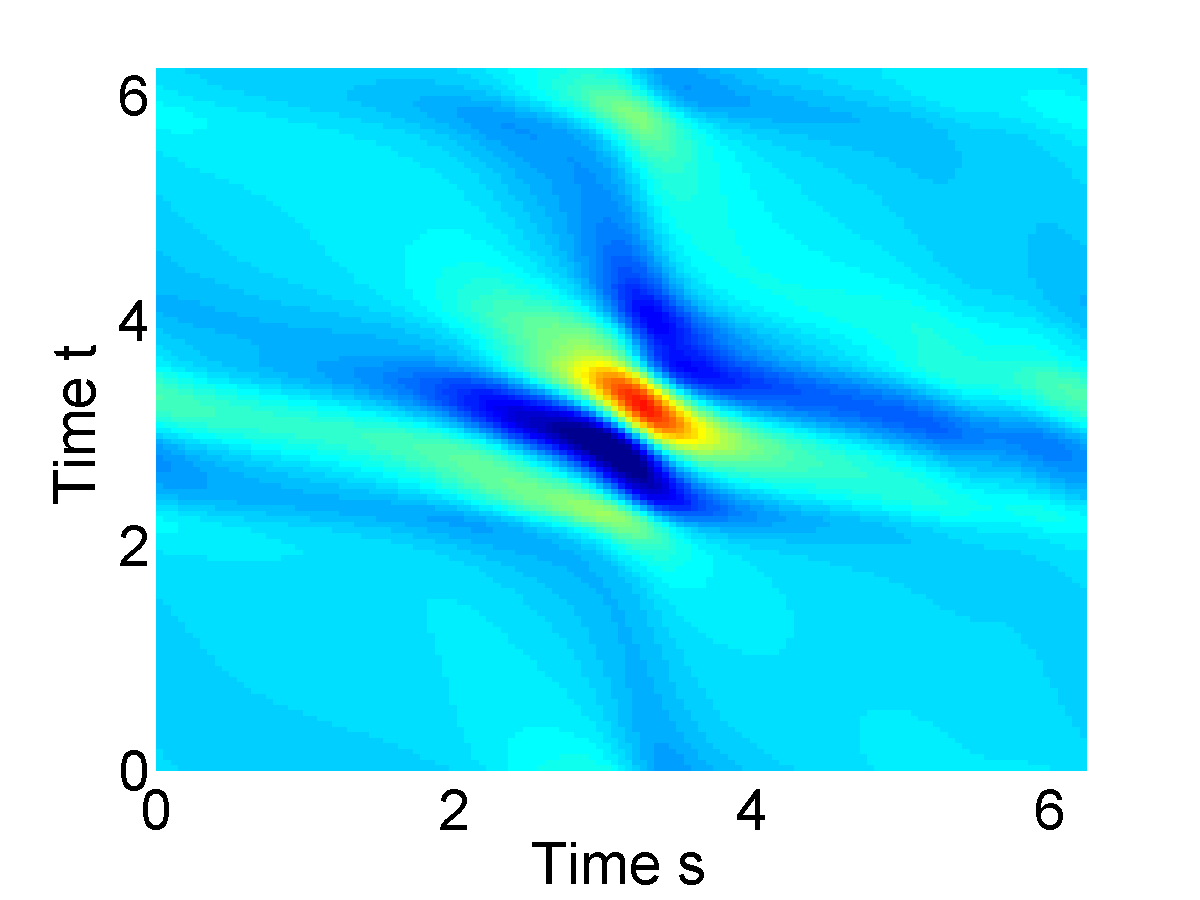}&\includegraphics[width=1.1in]{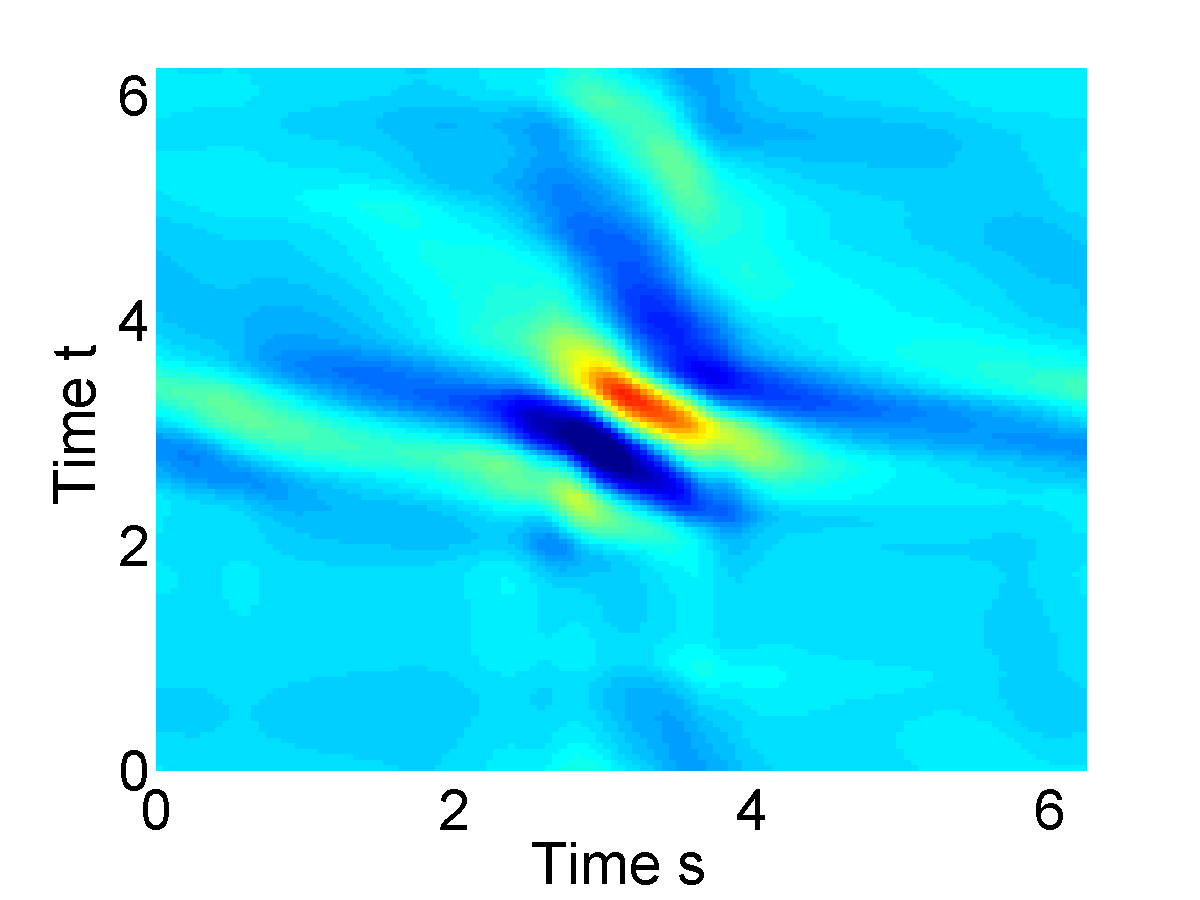}\\
\noalign{\vskip-10mm}
\hline
\noalign{\vskip 1mm}
\begin{tabular}{c}Remainder $r_{25}$\\\noalign{\vskip20mm}\end{tabular}&\includegraphics[width=1.1in]{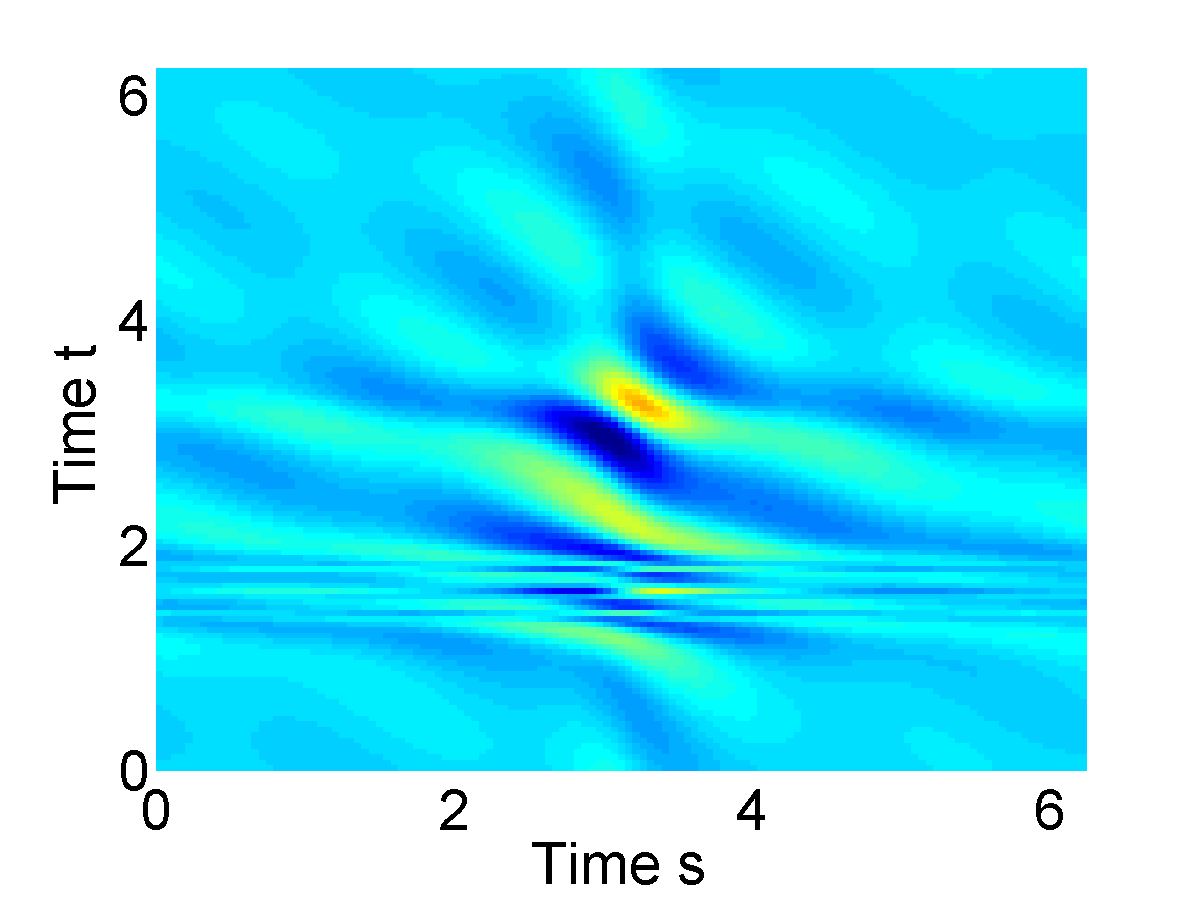}&\includegraphics[width=1.1in]{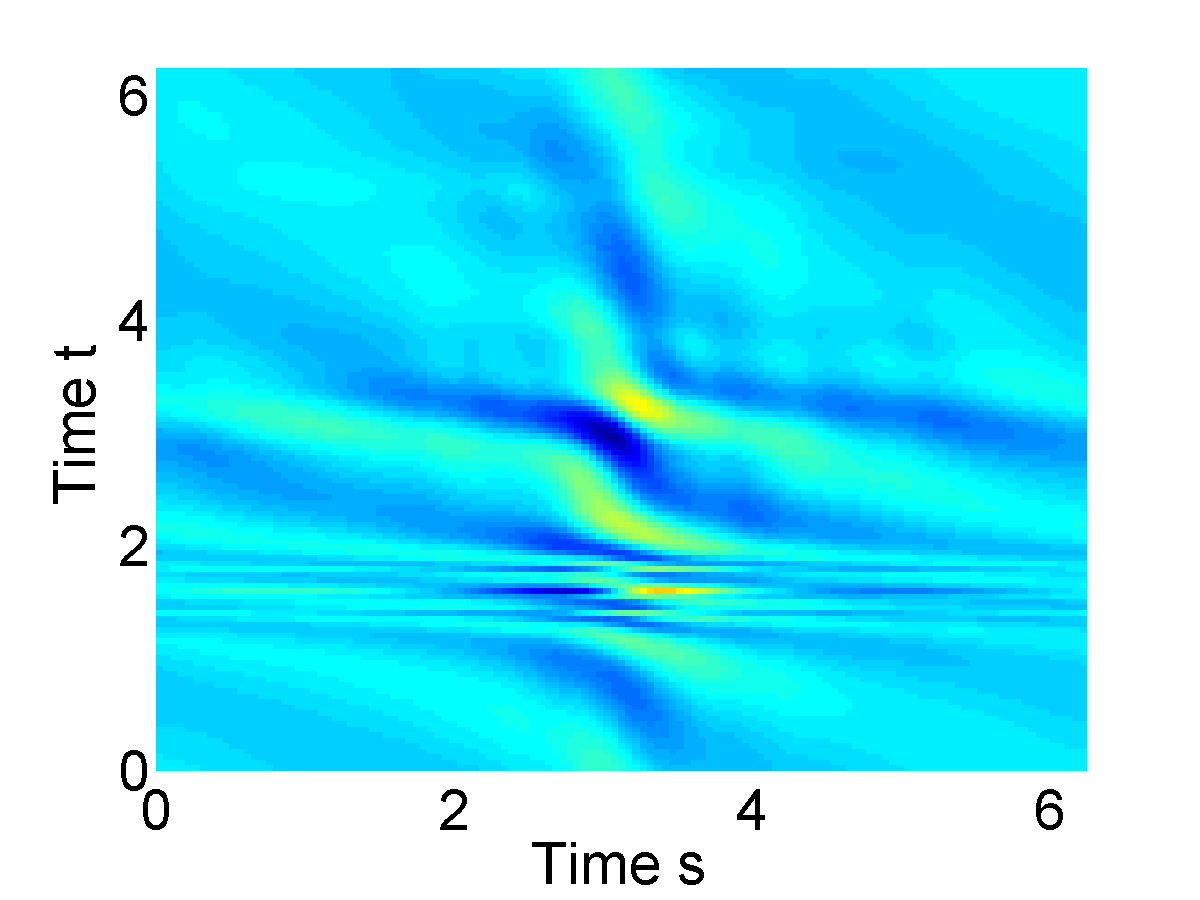}&\includegraphics[width=1.1in]{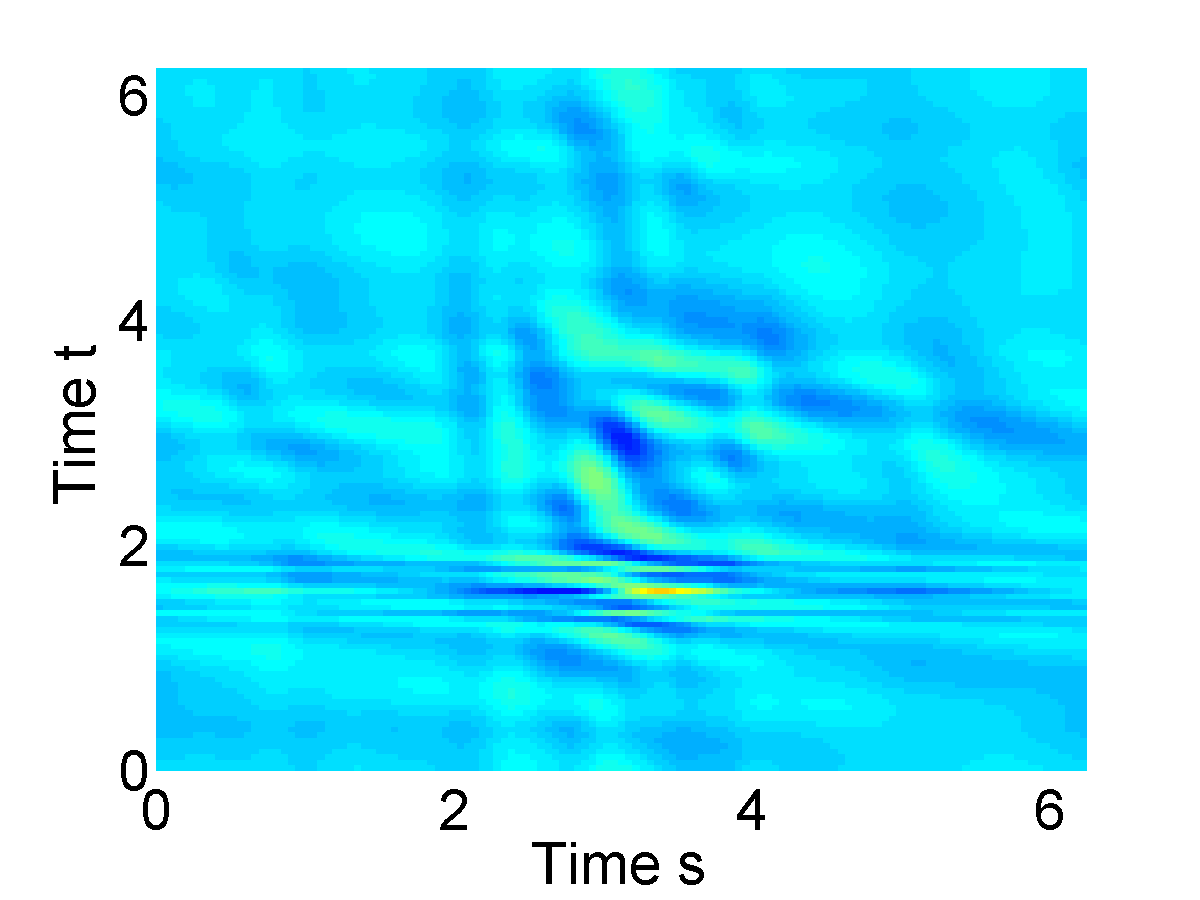}&\includegraphics[width=1.1in]{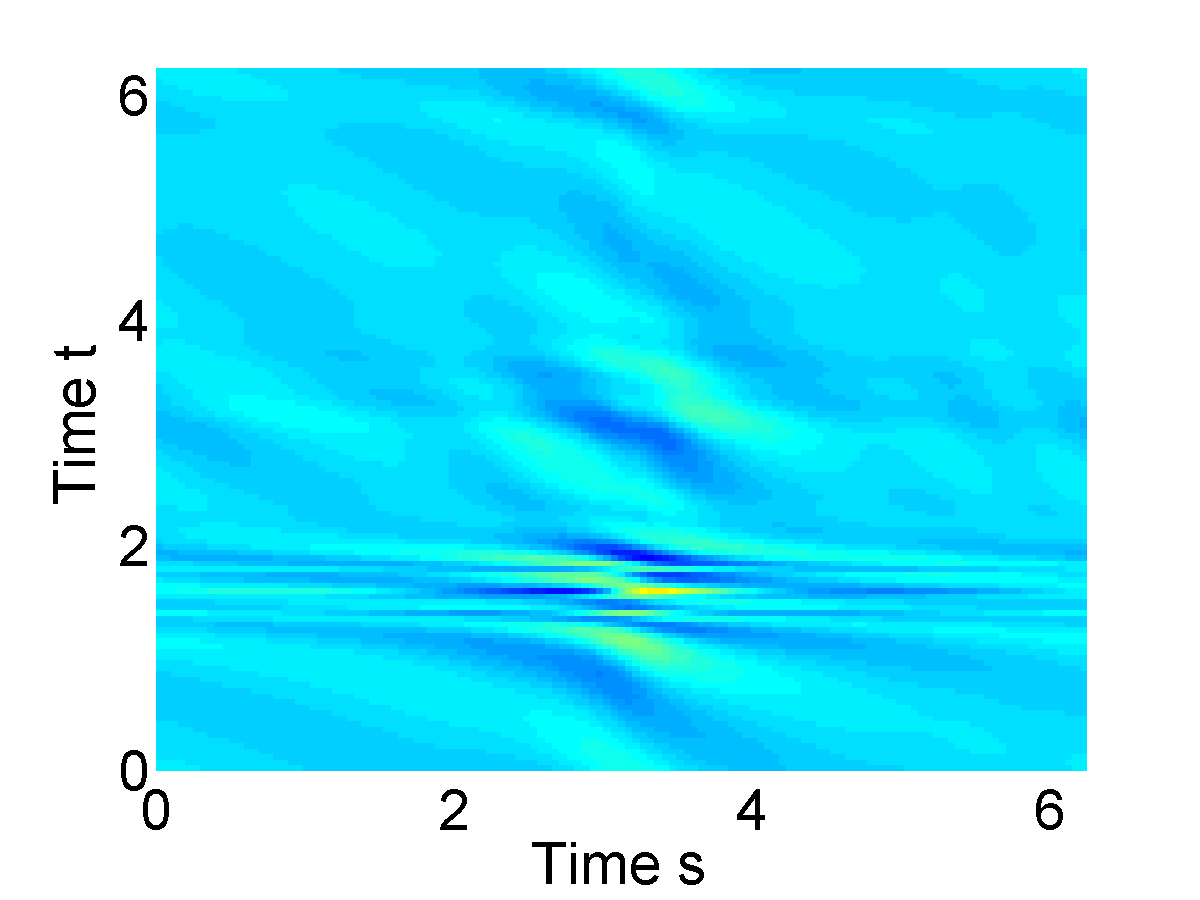}&\includegraphics[width=1.1in]{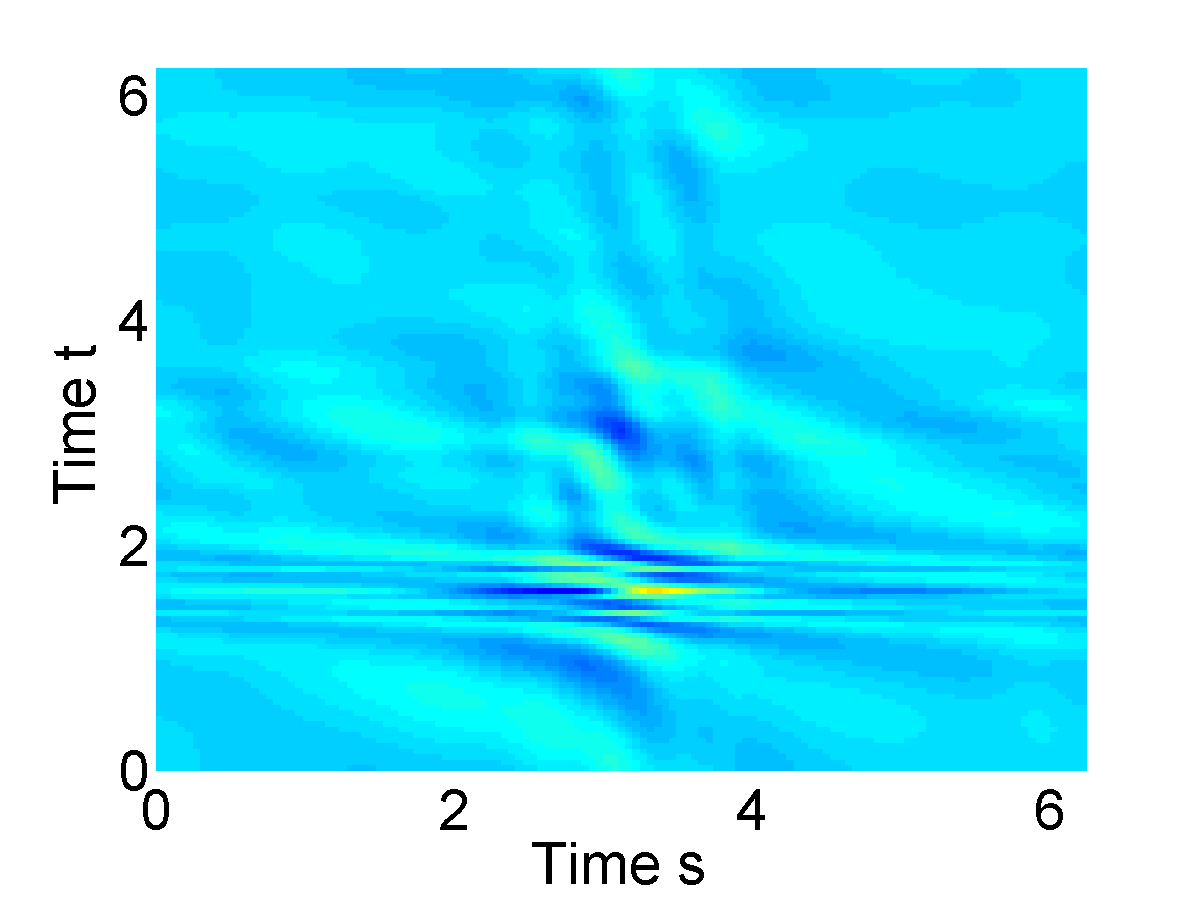}\\
\noalign{\vskip-10mm}
\hline
\end{tabular}
\end{table*}

For $f\in {\mathcal{H}}(\tilde{\mathcal{D}},M)$, it holds
\begin{equation}\begin{split}
|\langle f_{n},f\rangle |=&|\langle f_{n},\sum_{k=1}^\infty d_k\tilde B_k\rangle |\\\leqslant & \sum_{k=1}^\infty |d_k|\max_{k\geqslant 1} |\langle f_n,\tilde B_k\rangle|\\ \leqslant & \sum_{k=1}^\infty |d_k|\max_{\psi\in \mathcal{D}} |\langle f_n,B_n^{\psi}\rangle|\\=&N\max_{\psi\in \mathcal{D}} |\langle f_n,B_n^{\psi}\rangle|.\\
\end{split}\end{equation}
Namely,$$\max_{\psi\in \mathcal{D}} |\langle f_n,B_n^{\psi}\rangle|\geqslant \frac{|\langle f_{n},f\rangle |}{N}= \frac{\|f_{n}\|^2}{N}.$$

Then we have
\begin{equation}\begin{split}\|f_{n+1}\|^2\leqslant & \|f_{n}\|^2(1-\frac{\|f_{n}\|^2}{N^2}).\end{split}\end{equation}
Using Lemma \ref{amb} for $ t \equiv 1$, we obtain the newer error bound \begin{equation}\label{poganeb}\|f_n\|\leqslant \frac{N}{\sqrt{n}}.\end{equation}

\begin{remark}
Let $\mathcal{H}$ be a complex Hilbert space equipped with a dictionary $\mathcal{D}$. $\tilde{\mathcal{D}}$ is the complete dictionary induced from $\mathcal{D}$. Let $f\in \mathcal{H}(\mathcal{D},M)\cap\mathcal{H}(\mathcal{\tilde D},N)$. Let $M_f$, $N_f$ be respectively the infimum $M$ and the infimum $N$ such that $f\in \mathcal{H}(\mathcal{D},M)$ and $f\in \mathcal{H}(\mathcal{\tilde D},N)$. Then $N\leqslant M$. With regards to this, although (\ref{ogaeb}) and (\ref{poganeb}) are the same type of estimates, (\ref{poganeb}) is a better one.
\end{remark}
\section{Experiments on A Function With Two Variables}
To implement our algorithms, we take a signal of two variables as the tensor product of two 1D signals. Obviously, $H^2 (\mathbb{T})\otimes H^2 (\mathbb{T})\subset H^2 (\mathbb{T}^2)$. Aiming to take the singular measure into the consideration, we add a Dirac type point mass at $-i$ into a 1D signal. Then, the 2D function $f$ is given as
\begin{equation}f(z,w)=f_1\otimes f_2, \quad z=e^{it},w=e^{is},\end{equation}
where
\begin{equation}
f_1(z)=\frac{4z^2(1+0.02z)}{(1+0.7z)}e^{\frac{z+i}{z-i}},\quad
f_2(w)=\frac{4w^2(1+0.02w)}{(1+0.7w)}.
\end{equation}
\begin{figure}[h]
\centering
\includegraphics[width=1.1in]{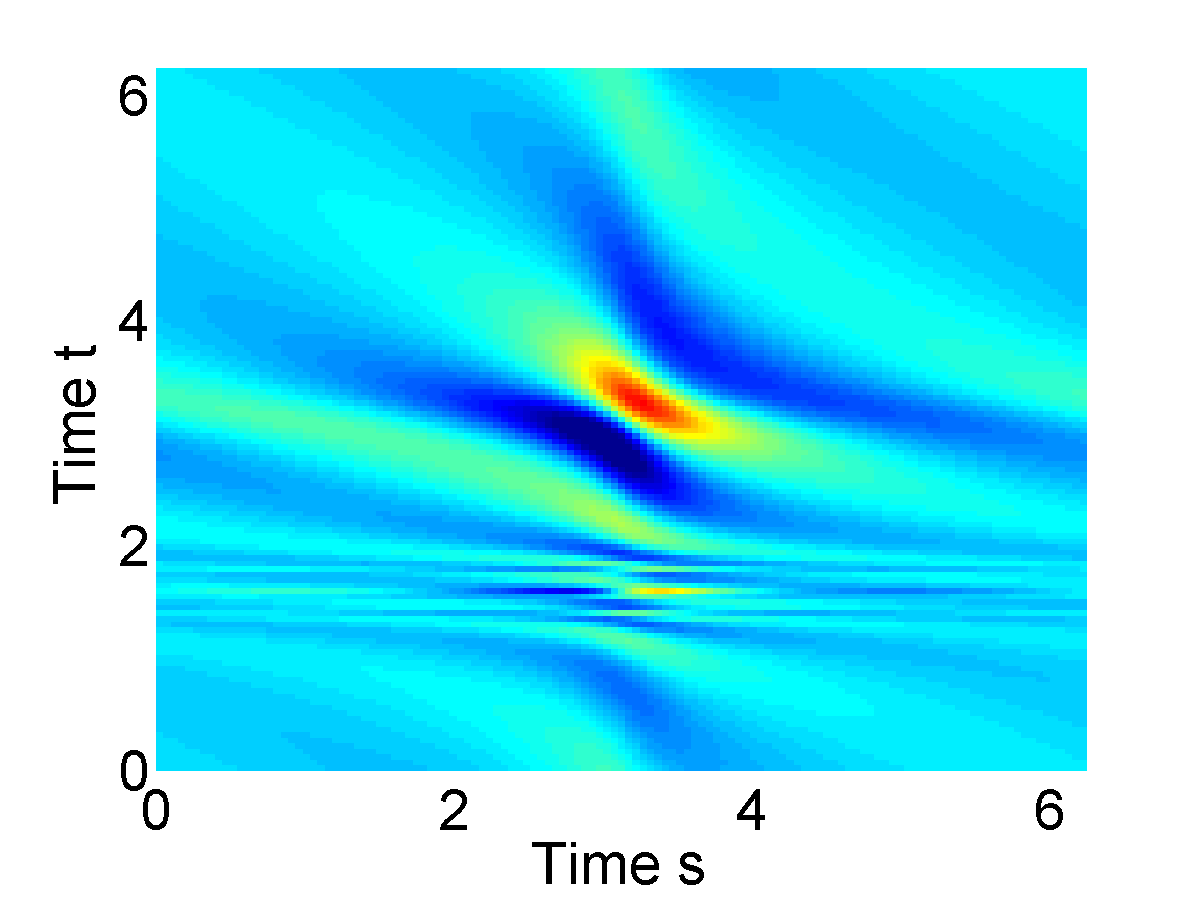}\hspace{1cm}
\includegraphics[width=1.1in]{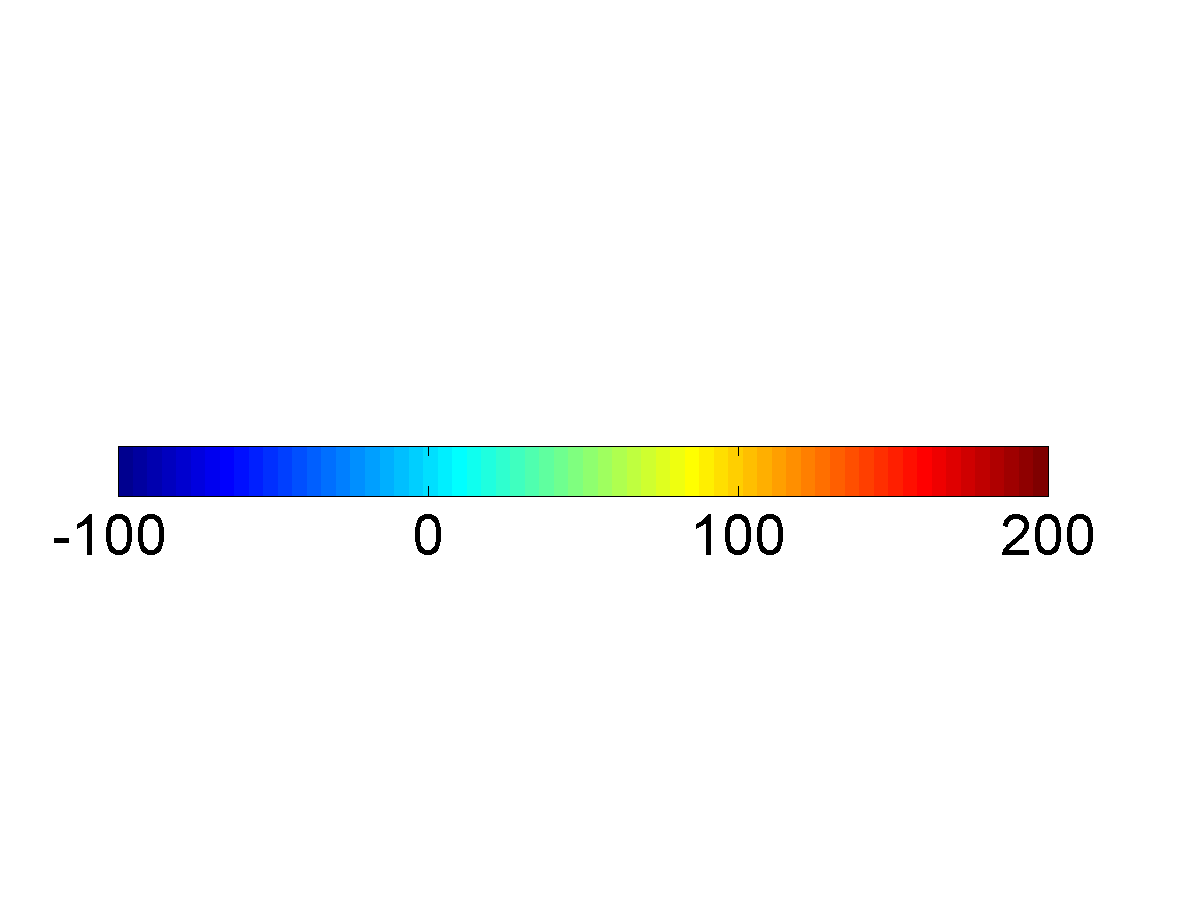}
\caption{The real part of the original function with corresponding colormap}
\label{f5}
\end{figure}
\begin{figure}[t]
\centering
\includegraphics[width=2.9in]{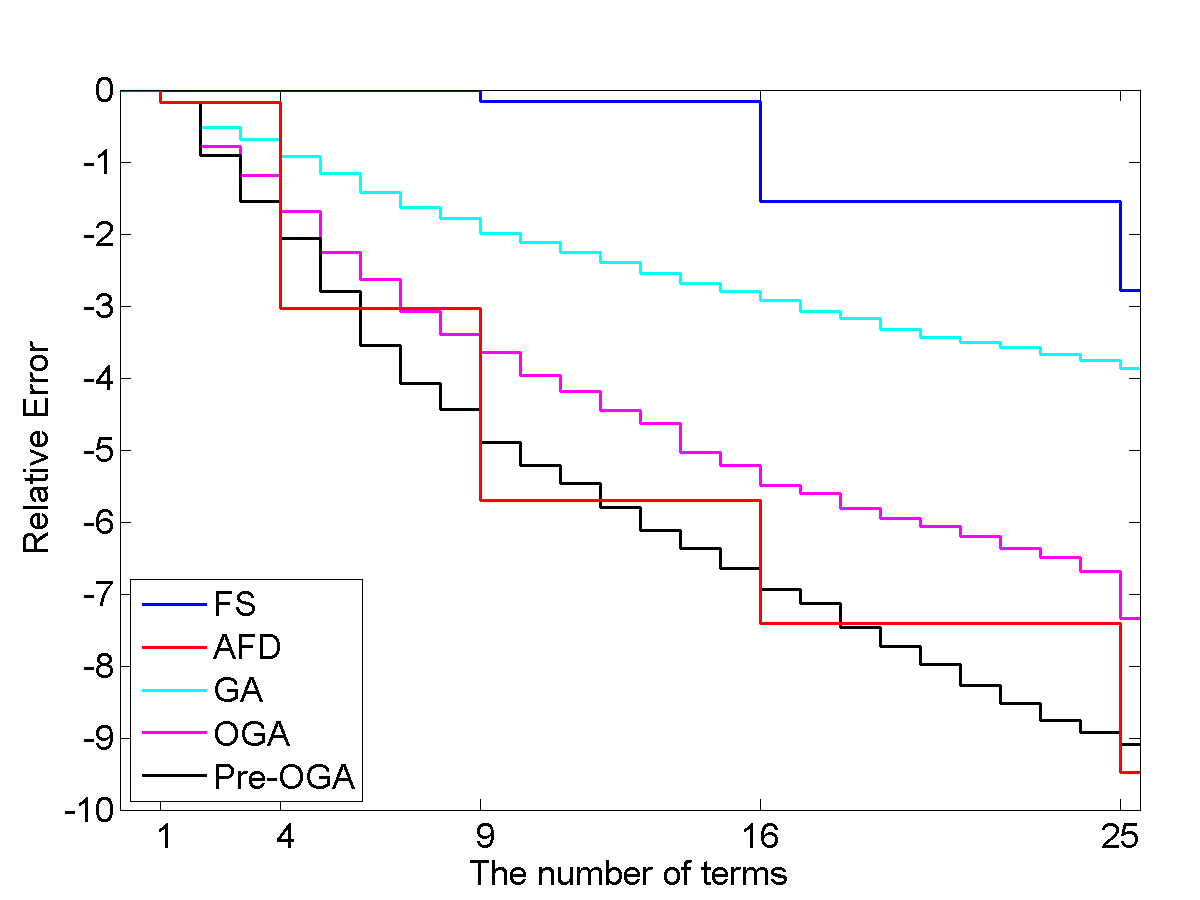}
\caption{Relative Error $r_N$}
\label{re5}
\end{figure}

$f (z, w) $ is sampled by $128\times 128$ points on $ \mathbb{T}^2$ evenly. The real part of the discrete signal $\bold F$ is shown in Fig \ref{f5}. 
By utilizing 193 grid points spaced evenly in the unit disc as options of the parameters, the $N$-term approximation result $\bold F_N$ from five algorithms be computed numerically in Matlab2012b. All results take the same colorbar with the original signal. Note that within these decompositions, at $m$-th step Product AFD yields $m^2$ terms that other algorithms need $m^2$ steps to reach. FD is implemented by setting all parameters in Algorithm \ref{2dafd} to zeros. Table \ref{f5img} shows the real parts of $\bold F_1,\;\bold F_9$ and $\bold F_{25}$ from five algorithms respectively. Remainders after 25 steps are listed at the bottom of Table \ref{f5img}. dB is used as the units of the relative error $r_N$.

From $r_N$ in Figure \ref{re5}, we find that $r_N$ of FD begins to fall from $9$-terms. This can also be observed from the expression of $f(z, w)$. Four adaptive algorithms select the same dictionary element at the first step, which is better than the one from FD. It accounts for the selection of the basic functions in $\mathcal{D}_S$ which facilitates the acquisition of proper energy. The approximation of Product AFD surpasses all the other algorithms. Among all algorithms, GA performs at an average level. It can be said that Product AFD is the best one among these methods from the standpoint of energy approximation. This may also be caused by the tensor product structure of the signal $f(z, w)$.
\section{Experiments on An Image}

We conduct experiments of image signal decompositions with Product AFD, 2D-Pre-OGA, GA, OGA, and FD. Through visual perception and three numerical indicators, the reconstruction capacity is evaluated subjectively and objectively. 

Histogram comparison is conducted among the normalized histogram of the approximation results. The difference between the histogram of the reconstructed image $H_1$ and one of the original image $H_2$ can be quantitated by Bhattacharyya distance $d$ \cite{Bhattacharyya1943,Kailath1967},
where $$d=\sqrt{1-\frac{\sum\sqrt{H_1H_2}}{\sqrt{\sum H_1\sum H_2}}}.$$

A histogram views image pixel as being independent of each other since it omits the corresponding relationship. Such objective evaluation then appears to have great difference from subjective evaluation \cite{Gonzalez2009, Eskicioglu1995}. To compensate this situation, two more parameters are used to evaluate the quality of reconstructed image: Peak Signal to Noise Rate (PSNR) and Quality Assessment (QA, viz. Mean Structural SIMilarity (MSSIM)). The PSNR views image pixels as independent variables, too. QA, however, offers a global measurement. 

The definition of PSNR is
$$PSNR=10\times \log\frac{L\times L}{MSE},$$
$$MSE=\frac{1}{M\times N}\sum_{i=1}^N\sum_{j=1}^M(\bold F[i,j]-\tilde{\bold F}_N[i,j])^2.$$ $(M,N)$ is the size of image $\bold F$. $L$ is the maximum value of the image pixels. Because the original image 'bird' is 8 bits image, $L= 255$. The definition of QA can be found in \cite{Wang2004}. We take the weight in QA an $11\times 11$ circular-symmetric Gaussian weighting function with $\sum_{i=1}^N\omega_i=1$.

\begin{table}[t]
\caption{Bhattacharyya distance between the histograms of the $N$-terms approximation results and the original image}
\label{hist}
\centering
\def\arraystretch{1.5}
\begin{tabular}{|c|c|c|c|c|c|}
\hline
$N$&GA&FD&Product AFD&OGA&2D-Pre-OGA\\
\hline
$1$ & 0.6379 & 0.6213 & 0.6436 &0.6376 & 0.6376\\
\hline 
$16$& 0.7637 & 0.7696 & 0.7901 & 0.8107 & 0.7992\\
\hline
$64$& 0.8178& 0.8194& 0.8247&0.8374 & 0.8456\\
\hline
$256$& 0.8368 & 0.8516 & 0.8678 & 0.8684 & 0.8763\\
\hline
\end{tabular}
\end{table}

\begin{figure*}[!th]
\centering
\begin{minipage}[b]{.15\textwidth}
\includegraphics[width=1in]{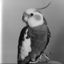}
\vspace{1cm}
\caption{The original image: bird}
\label{bird}
\end{minipage}
\begin{minipage}[b]{.4\textwidth}
\centering
\includegraphics[width=2.9in]{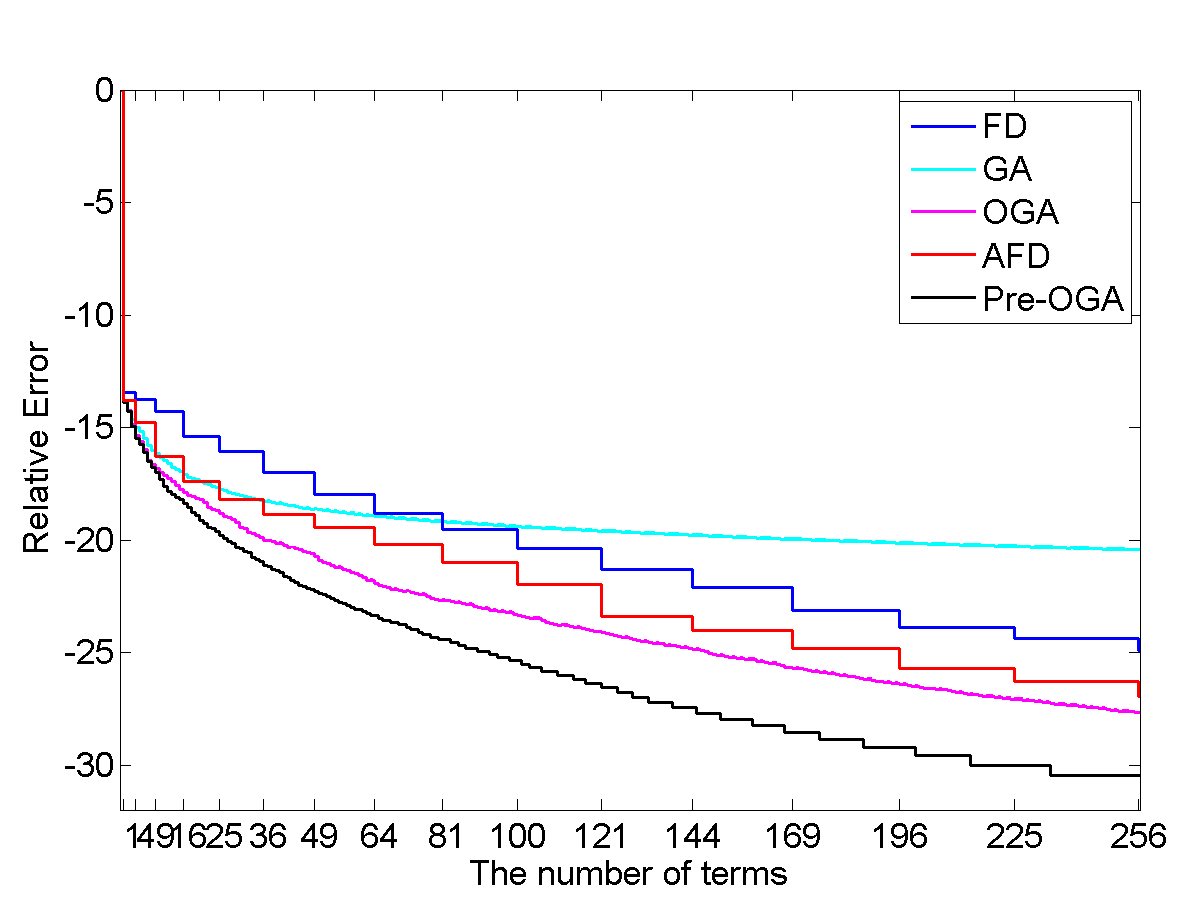}
\caption{PSNR}
\label{re6}
\end{minipage}
\begin{minipage}[b]{.4\textwidth}
\centering
\includegraphics[width=2.9in]{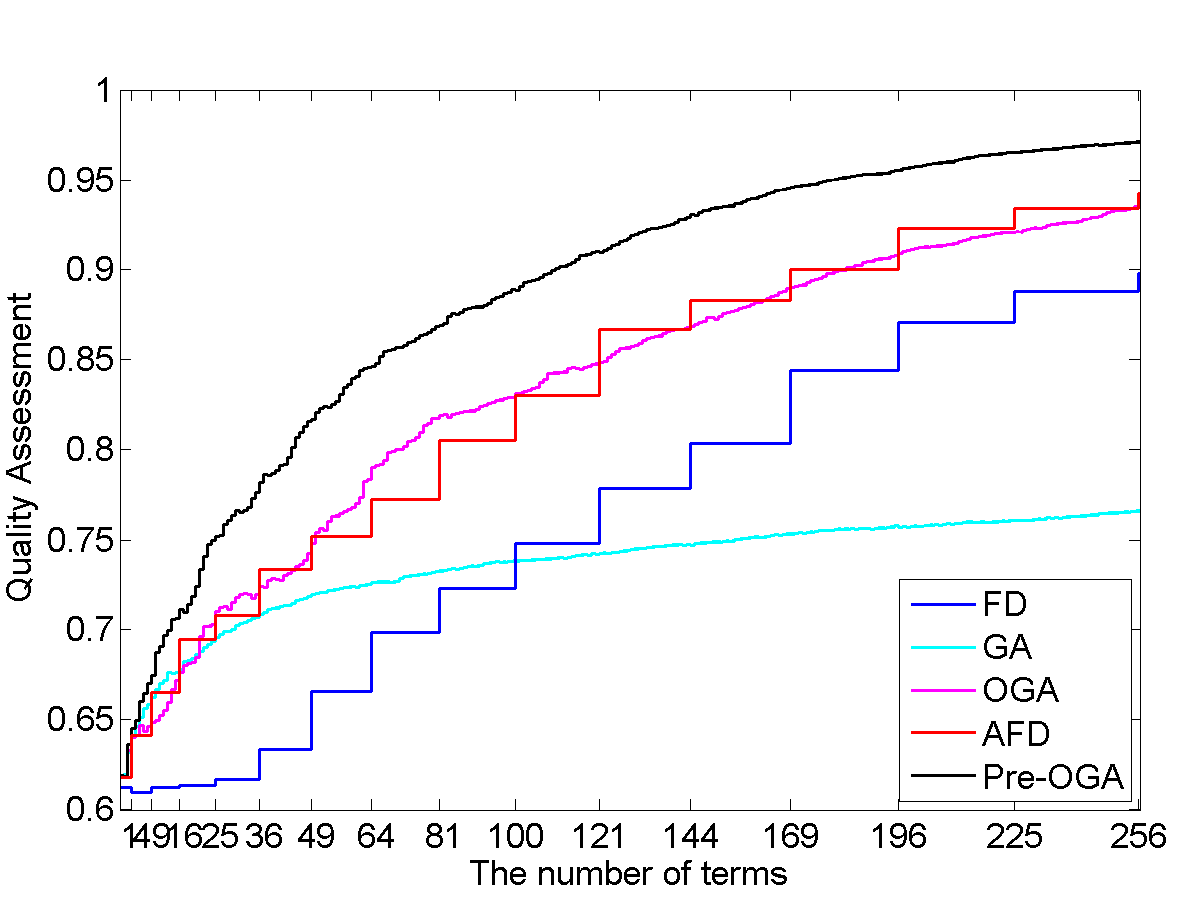}
\caption{QA}
\label{qa6}
\end{minipage}
\end{figure*}

The original image is a $64\times64$ PNG photo as Fig. \ref{bird}.
The reconstruction results $\tilde{\bold F}_1,\;\tilde{\bold F}_{16},\;\tilde{\bold F}_{64}$ and $\tilde{\bold F}_{256}$ are shown in Table \ref{image}. The remainders after 256 steps from five methods are given by its logarithm.

\begin{table*}[t]
\centering
\caption{$\bold F_1,\;\bold F_{16},\;\bold F_{64}$, $\bold F_{256}$, and the logarithm of the remainder after 256-terms approximation}
\label{image}
\begin{tabular}{cccccc}
\hline
Algorithm&FS&GA&OGA&Product AFD&Pre-OGA\\
\hline
\noalign{\vskip 1mm}
\begin{tabular}{c}$\tilde{\bold F}_1$\\\noalign{\vskip20mm}\end{tabular} &\includegraphics[width=1in]{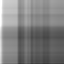}&\includegraphics[width=1in]{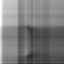}&\includegraphics[width=1in]{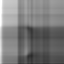}&\includegraphics[width=1in]{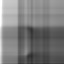}&\includegraphics[width=1in]{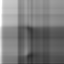}\\
\noalign{\vskip-10mm}
\hline
\noalign{\vskip 1mm}
\begin{tabular}{c}$\tilde{\bold F}_{16}$\\\noalign{\vskip20mm}\end{tabular}&\includegraphics[width=1in]{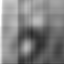}&\includegraphics[width=1in]{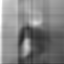}&\includegraphics[width=1in]{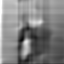}&\includegraphics[width=1in]{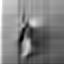}&\includegraphics[width=1in]{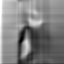}\\
\noalign{\vskip-10mm}
\hline
\noalign{\vskip 1mm}
\begin{tabular}{c}$\tilde{\bold F}_{64}$\\\noalign{\vskip20mm}\end{tabular}&\includegraphics[width=1in]{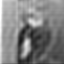}&\includegraphics[width=1in]{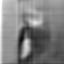}&\includegraphics[width=1in]{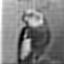}&\includegraphics[width=1in]{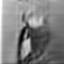}&\includegraphics[width=1in]{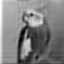}\\
\noalign{\vskip-10mm}
\hline
\noalign{\vskip 1mm}
\begin{tabular}{c}$\tilde{\bold F}_{256}$\\\noalign{\vskip20mm}\end{tabular}&\includegraphics[width=1in]{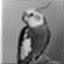}&\includegraphics[width=1in]{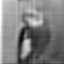}&\includegraphics[width=1in]{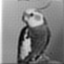}&\includegraphics[width=1in]{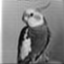}&\includegraphics[width=1in]{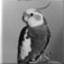}\\
\noalign{\vskip-10mm}
\hline\noalign{\vskip 1mm}
\begin{tabular}{c} $log_{10}(r_{256})$\\\noalign{\vskip18mm}\end{tabular}&\includegraphics[width=1in]{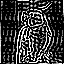}&\includegraphics[width=1in]{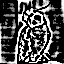}&\includegraphics[width=1in]{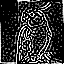}&\includegraphics[width=1in]{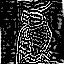}&\includegraphics[width=1in]{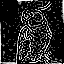}\\
\noalign{\vskip-10mm}
\hline
\end{tabular}
\end{table*}

From the direct observation of these approximation results, the first term approximations of these algorithms coincide with each other, but except FD. The stripes on $\tilde{\bold F}_1$ are generated from $\bold F_0$ and $\bold G_0$ in equation (\ref{Fc0}). The approximation results of 2D-Pre-OGA become more and more clear with increasing of $N$. GA generates blur results even with $256$-terms approximation. The Bhattacharyya distances between the histograms are shown in Table \ref{hist}. Histogram comparison is a convenient method of easy computation. PSNR and QA are shown in Fig. \ref{re6} and \ref{qa6}, respectively. After the approximations with 256 components, Pre-OGA reaches 30dB PSNR as the best result. Similarly, QA of 2D-Pre-OGA maintains its advantage. Product AFD and OGA are very close referring to QA.

\section{Conclusion}
This paper includes mainly three correlative works of multi-dimensional AFDs. In the first part, we reformulate the AFD theory in view of their algorithms, on the 2-torus. We propose the computerized numerical realization of Product AFD and 2D-Pre-OGA. The numerical realization of Product AFD and 2D-Pre-OGA are phrased as, respectively, Algorithms \ref{2dafd} and Algorithm \ref{poga}. Algorithm \ref{realsig} is preserved to the realization of applying the mentioned algorithms to real-valued signals. The idea of these numerical methods can be generalized to the numerical computation of multi-dimensional AFDs. In the second part, Product AFD, 2D-Pre-OGA, GA, OGA, and FD have been illustrated by toy and real image examples. The experiment results show that for a function with singularity and decentralized energy distribution, FD cannot yield a satisfactory convergence. It is possible that FD produces better results than that of the GA with sufficient terms. The orthogonalization involving dictionary elements selections, that is OGA, Product AFD, and 2D-Pre-OGA achieve better approximation in adaptive approximations. The comparisons drawn from the experiments indicate that Product AFD and Pre-OGA are preferable approximation methods. Without a Gram-Schmidt orthonormalization progress, Product AFD decomposes a signal into a rational orthogonal system, which can be implemented from an explicit expression with less numerical instability and less computation complexity. Intrinsically, Product AFD inherited non-negative frequency processing and fast energy convergence rate from 1D-AFD. In the experiments, the 2D-Pre-OGA generates better results, merely at a cost of time. The last part but not the least, we proposed new and sharper estimations for the error bounds of OGA and Pre-OGA. Through comparative analysis, the weak type of Pre-OGA can reach $ t =1$ during approximation, but the weak type of GA and OGA cannot.


%

%

\section*{Acknowledgment}
The first author would like to thank Xing-min Li, Li-ming Zhang, Wei-Xiong Mai, Wei Wu, Yang Wang for helpful discussions.

\ifCLASSOPTIONcaptionsoff
\newpage
\fi

\bibliographystyle{IEEEtran}
\bibliography{ref}

\begin{thebibliography}{10}
\providecommand{\url}[1]{#1}
\csname url@samestyle\endcsname
\providecommand{\newblock}{\relax}
\providecommand{\bibinfo}[2]{#2}
\providecommand{\BIBentrySTDinterwordspacing}{\spaceskip=0pt\relax}
\providecommand{\BIBentryALTinterwordstretchfactor}{4}
\providecommand{\BIBentryALTinterwordspacing}{\spaceskip=\fontdimen2\font plus
\BIBentryALTinterwordstretchfactor\fontdimen3\font minus
  \fontdimen4\font\relax}
\providecommand{\BIBforeignlanguage}[2]{{%
\expandafter\ifx\csname l@#1\endcsname\relax
\typeout{** WARNING: IEEEtran.bst: No hyphenation pattern has been}%
\typeout{** loaded for the language `#1'. Using the pattern for}%
\typeout{** the default language instead.}%
\else
\language=\csname l@#1\endcsname
\fi
#2}}
\providecommand{\BIBdecl}{\relax}
\BIBdecl

\bibitem{Qian2011}
T.~Qian and Y.-B. Wang, ``Adaptive fourier series—a variation of greedy
  algorithm,'' \emph{Advances in Computational Mathematics}, vol.~34, no.~3,
  pp. 279--293, 2011.

\bibitem{Qian2010}
T.~Qian, ``Intrinsic mono-component decomposition of functions: An advance of
  fourier theory,'' \emph{Mathematical Methods in the Applied Sciences},
  vol.~33, no.~7, pp. 880--891, 2010.

\bibitem{Qian2016}
------, ``Two-dimensional adaptive fourier decomposition,'' \emph{Mathematical
  Methods in the Applied Sciences}, vol.~39, no.~10, pp. 2431--2448, 2016.

\bibitem{Qian2014}
\BIBentryALTinterwordspacing
------, ``Cyclic afd algorithm for the best rational approximation,''
  \emph{Mathematical Methods in the Applied Sciences}, vol.~37, no.~6, pp.
  846--859, 2014. [Online]. Available: \url{http://dx.doi.org/10.1002/mma.2843}
\BIBentrySTDinterwordspacing

\bibitem{Mallat1993}
S.~G. Mallat and Z.~Zhang, ``Matching pursuits with time-frequency
  dictionaries,'' \emph{IEEE Transactions on Signal Processing}, vol.~41,
  no.~12, pp. 3397--3415, 1993.

\bibitem{Davis1994}
G.~Davis, ``Adaptive nonlinear approximations,'' Ph.D. dissertation, New York
  University, 1994.

\bibitem{DeVore1996}
R.~A. DeVore and V.~N. Temlyakov, ``Some remarks on greedy algorithms,''
  \emph{Advances in Computational Mathematics}, vol.~5, no.~1, pp. 173--187,
  1996.

\bibitem{Temlyakov2000}
\BIBentryALTinterwordspacing
V.~Temlyakov, ``Weak greedy algorithms,'' \emph{Advances in Computational
  Mathematics}, vol.~12, no. 2-3, pp. 213--227, 2000. [Online]. Available:
  \url{http://dx.doi.org/10.1023/A3A1018917218956}
\BIBentrySTDinterwordspacing

\bibitem{Tem1}
------, \emph{Greedy approximation}.\hskip 1em plus 0.5em minus 0.4em\relax
  Cambridge University Press, 2011.

\bibitem{Mi2012}
W.~Mi and T.~Qian, ``Frequency-domain identification: An algorithm based on an
  adaptive rational orthogonal system,'' \emph{Automatica}, vol.~48, no.~6, pp.
  1154--1162, 2012.

\bibitem{Qian2013}
\BIBentryALTinterwordspacing
T.~Qian and Y.~Wang, ``Remarks on adaptive fourier decomposition,''
  \emph{International Journal of Wavelets, Multiresolution and Information
  Processing}, vol.~11, no.~01, p. 1350007, 2013. [Online]. Available:
  \url{http://www.worldscientific.com/doi/abs/10.1142/S0219691313500070}
\BIBentrySTDinterwordspacing

\bibitem{Mi2016}
W.~Mi, T.~Qian, and S.~Li, ``Basis pursuit for frequency-domain
  identification,'' \emph{Mathematical Methods in the Applied Sciences},
  vol.~39, no.~3, pp. 498--507, 2016.

\bibitem{Mo2014}
Y.~Mo and T.~Qian, ``Support vector machine adapted tikhonov regularization
  method to solve dirichlet problem,'' \emph{Applied Mathematics and
  Computation}, vol. 245, pp. 509--519, 2014.

\bibitem{Mi2014}
W.~Mi and T.~Qian, ``On backward shift algorithm for estimating poles of
  systems,'' \emph{Automatica}, vol.~50, no.~6, pp. 1603 -- 1610, 2014.

\bibitem{Mo2015}
Y.~Mo, T.~Qian, and W.~Mi, ``Sparse representation in szegő kernels through
  reproducing kernel hilbert space theory with applications,''
  \emph{International Journal of Wavelets, Multiresolution and Information
  Processing}, vol.~13, no.~04, p. 1550030, 2015.

\bibitem{Takenaka1925}
S.~Takenaka, ``On the orthogonal functions and a new formula of
  interpolation,'' in \emph{Japanese Journal of Mathematics: Transactions and
  Abstracts}, vol.~2.\hskip 1em plus 0.5em minus 0.4em\relax The Mathematical
  Society of Japan, 1925, pp. 129--145.

\bibitem{Walsh1935}
J.~L. Walsh, \emph{Interpolation and approximation by rational functions in the
  complex domain}.\hskip 1em plus 0.5em minus 0.4em\relax American Mathematical
  Soc., 1935, vol.~20.

\bibitem{Bultheel1999}
A.~Bultheel, \emph{Orthogonal rational functions}.\hskip 1em plus 0.5em minus
  0.4em\relax Cambridge University Press, 1999.

\bibitem{Akccay2001}
H.~Ak{\c{c}}ay, ``On the uniform approximation of discrete-time systems by
  generalized fourier series,'' \emph{IEEE Transactions on Signal Processing},
  vol.~49, no.~7, pp. 1461--1467, 2001.

\bibitem{heuberger2005}
P.~S. Heuberger, P.~M. van~den Hof, and B.~Wahlberg, \emph{Modelling and
  identification with rational orthogonal basis functions}.\hskip 1em plus
  0.5em minus 0.4em\relax Springer Science \& Business Media, 2005.

\bibitem{Szasz1953}
\BIBentryALTinterwordspacing
O.~Sz{\'a}sz, ``On closed sets of rational functions,'' \emph{Annali di
  Matematica Pura ed Applicata}, vol.~34, no.~1, pp. 195--218, 1953. [Online].
  Available: \url{http://dx.doi.org/10.1007/BF02415331}
\BIBentrySTDinterwordspacing

\bibitem{szeg1939}
G.~Szeg\"{o}, \emph{Orthogonal polynomials}.\hskip 1em plus 0.5em minus
  0.4em\relax American Mathematical Soc., 1939, vol.~23.

\bibitem{lee1961}
Y.~W. Lee, ``Statistical theory of communication,'' \emph{American Journal of
  Physics}, vol.~29, no.~4, pp. 276--278, 1961.

\bibitem{Broome1965}
\BIBentryALTinterwordspacing
P.~W. Broome, ``Discrete orthonormal sequences,'' \emph{Journal of the ACM},
  vol.~12, no.~2, pp. 151--168, Apr. 1965. [Online]. Available:
  \url{http://doi.acm.org/10.1145/321264.321265}
\BIBentrySTDinterwordspacing

\bibitem{Cohen1995}
L.~Cohen, \emph{Time-frequency analysis}.\hskip 1em plus 0.5em minus
  0.4em\relax Prentice Hall PTR Englewood Cliffs, NJ:, 1995, vol. 778.

\bibitem{Weiss1962}
M.~Weiss and G.~Weiss, ``A derivation of the main results of the theory of hp
  spaces,'' \emph{Revista de la Uni{\'o}n Matem{\'a}tica Argentina}, vol.~20,
  pp. 63--71, 1962.

\bibitem{Nahon2000}
M.~R. Nahon, ``Phase evaluation and segmentation,'' Ph.D. dissertation, Yale
  University, New Haven, CT, USA, 2000.

\bibitem{Coifman2015a}
\BIBentryALTinterwordspacing
R.~R. Coifman and S.~Steinerberger, ``Nonlinear phase unwinding of functions,''
  \emph{Journal of Fourier Analysis and Applications}, p. 1–32, 2016.
  [Online]. Available: \url{http://dx.doi.org/10.1007/s00041-016-9489-3}
\BIBentrySTDinterwordspacing

\bibitem{Qian2012}
T.~Qian, W.~Sprößig, and J.~Wang, ``Adaptive fourier decomposition of
  functions in quaternionic hardy spaces,'' \emph{Mathematical Methods in the
  Applied Sciences}, vol.~35, no.~1, pp. 43--64, 2012.

\bibitem{Qian2014c}
T.~Qian, J.~Wang, and Y.~Yang, ``Matching pursuits among shifted cauchy kernels
  in higher-dimensional spaces,'' \emph{Acta Mathematica Scientia}, vol.~34,
  no.~3, pp. 660 -- 672, 2014.

\bibitem{Alpay2016}
\BIBentryALTinterwordspacing
D.~Alpay, F.~Colombo, T.~Qian, and I.~Sabadini, ``Adaptative decomposition: The
  case of the {D}rury--{A}rveson space,'' \emph{Journal of Fourier Analysis and
  Applications}, Oct 2016. [Online]. Available:
  \url{http://dx.doi.org/10.1007/s00041-016-9508-4}
\BIBentrySTDinterwordspacing

\bibitem{Alpay2017}
------, ``Adaptive orthonormal systems for matrix-valued functions,''
  \emph{Proceedings of the American Mathematical Society}, vol. 145, pp.
  2089--2106, 2017.

\bibitem{Rudin1969}
\BIBentryALTinterwordspacing
W.~Rudin, \emph{Function theory in polydiscs}, ser. Mathematics lecture note
  series.\hskip 1em plus 0.5em minus 0.4em\relax W. A. Benjamin, 1969, no.~41.
  [Online]. Available: \url{https://books.google.com/books?id=9waoAAAAIAAJ}
\BIBentrySTDinterwordspacing

\bibitem{Bracewell1986}
R.~N. Bracewell, \emph{The Fourier transform and its applications}.\hskip 1em
  plus 0.5em minus 0.4em\relax McGraw-Hill New York, 1986, vol. 31999.

\bibitem{Hardy1915}
\BIBentryALTinterwordspacing
G.~H. Hardy, ``The mean value of the modulus of an analytic function,''
  \emph{Proceedings of the London Mathematical Society}, vol. s2$\_$14, no.~1,
  pp. 269--277, 1915. [Online]. Available:
  \url{http://dx.doi.org/10.1112/plms/s2$\_$14.1.269}
\BIBentrySTDinterwordspacing

\bibitem{Garnett1981}
J.~B. Garnett, \emph{Bounded analytic functions}.\hskip 1em plus 0.5em minus
  0.4em\relax Academic press, 1981, vol.~96.

\bibitem{Zygmund2002}
A.~Zygmund, \emph{Trigonometric series}.\hskip 1em plus 0.5em minus 0.4em\relax
  Cambridge university press, 2002, vol.~1.

\bibitem{DeVore1998}
R.~A. DeVore, ``Nonlinear approximation,'' \emph{Acta Numerica}, vol.~7, p.
  51–150, 1998.

\bibitem{Aronszajn1950}
N.~Aronszajn, ``Theory of reproducing kernels,'' \emph{Transactions of the
  American mathematical society}, vol.~68, no.~3, pp. 337--404, 1950.

\bibitem{Bhattacharyya1943}
\BIBentryALTinterwordspacing
A.~Bhattacharyya, ``On a measure of divergence between two statistical
  populations defined by their probability distributions,'' \emph{Bulletin of
  Calcutta Mathematical Society}, vol.~35, no.~1, pp. 99--109, 1943. [Online].
  Available: \url{http://ci.nii.ac.jp/naid/10027606363/en/}
\BIBentrySTDinterwordspacing

\bibitem{Kailath1967}
T.~Kailath, ``The divergence and bhattacharyya distance measures in signal
  selection,'' \emph{IEEE Transactions on Communication Technology}, vol.~15,
  no.~1, pp. 52--60, February 1967.

\bibitem{Gonzalez2009}
R.~C. Gonzalez and R.~E. Woods, \emph{Digital Image Processing}.\hskip 1em plus
  0.5em minus 0.4em\relax Pearson Education, 2009.

\bibitem{Eskicioglu1995}
A.~M. Eskicioglu and P.~S. Fisher, ``Image quality measures and their
  performance,'' \emph{IEEE Transactions on Communications}, vol.~43, no.~12,
  pp. 2959--2965, Dec 1995.

\bibitem{Wang2004}
Z.~Wang, A.~C. Bovik, H.~R. Sheikh, and E.~P. Simoncelli, ``Image quality
  assessment: from error visibility to structural similarity,'' \emph{IEEE
  Transactions on Image Processing}, vol.~13, no.~4, pp. 600--612, 2004.

\end{thebibliography}

\end{document}